\documentclass{amsart}
\usepackage[latin1]{inputenc}
\makeatletter

\theoremstyle{plain}    
\newtheorem{thm}{Theorem}[section]
\numberwithin{equation}{section} 
\numberwithin{figure}{section} 
\theoremstyle{plain}    
\theoremstyle{plain}    
\newtheorem{conj}[thm]{Conjecture} 
\newtheorem{lem}[thm]{Lemma} 
\theoremstyle{plain}    
\newtheorem{prop}[thm]{Proposition} 
\theoremstyle{plain}    
\newtheorem{Def}[thm]{Definition} 
\theoremstyle{remark}
\newtheorem{rem}[thm]{Remark}
\theoremstyle{remark}

\usepackage[arrow,matrix,curve,ps]{xy}
\xyoption{dvips}
\CompileMatrices

\makeatother

\begin{document}

\title[NUMERICALLY TRIVIAL FOLIATIONS AND NUMERICAL DIMENSION]{Numerically trivial foliations, Iitaka fibrations and the 
numerical dimension}


\author{Thomas Eckl}




\address{Thomas Eckl, Mathematisches Institut, Universität
  zu K\"oln, 50931 K\"oln, Germany}

\email{thomas.eckl@math.uni-koeln.de}

\urladdr{http://www.mi.uni-koeln.de/\~{}teckl}

\begin{abstract}
Modifying the notion of numerically trivial foliation of a pseudo-effective 
line bundle $L$ introduced by the author in \cite{Eck04b} (see also 
math.AG/0304312) it can be shown that the leaves of this foliation have 
codimension $\geq$ the numerical dimension of $L$ as defined by Boucksom, 
Demailly, Paun and Peternell, math.AG/0405285. Furthermore,
if the Kodaira dimension of $L$ equals its numerical dimension the 
Kodaira-Iitaka fibration is its numerically trivial foliation. Both statements
together yield a sufficient criterion for $L$ not being abundant.
\end{abstract}

\maketitle

\setcounter{section}{-1}

\section{Introduction}

\noindent In their seminal paper \cite{BDPP04} Boucksom, Demailly, Paun and 
Peternell introduced a numerical dimension for pseudo-effective 
$(1,1)$-classes on 
compact K\"ahler manifolds generalizing the numerical dimension of nef line 
bundles on projective 
manifolds. For this purpose they used Boucksom's moving intersection numbers 
\cite{Bou02} which can be defined as follows:
\begin{Def}
Let
$X$
be a compact K\"ahler manifold with K\"ahler form
$\omega$.
Let
$\alpha_1, \ldots ,\alpha_p \in H^{1,1}(X,\mathbb{R})$
be pseudo-effective classes and let
$\Theta$
be a closed positive current of bidimension
$(p,p)$.
Then the moving intersection number 
$(\alpha_1 \cdot \ldots \cdot \alpha_p \cdot \Theta)_{\geq 0}$
of the
$\alpha_i$
and
$\Theta$
is defined to be the limit when
$\epsilon > 0$
goes to
$0$
of 
\[ \sup \int_{X-F} (T_1+\epsilon\omega) \wedge \ldots \wedge 
   (T_p+\epsilon\omega) \wedge \Theta \]
where the 
$T_i$'s
run through all currents with analytic singularities in
$\alpha_i[-\epsilon\omega]$,
and 
$F$
is the union of the
$Sing(T_i)$.
\end{Def}

\noindent This may be used for
\begin{Def}
Let
$X$
be an $n$-dimensional compact K\"ahler manifold. Then the numerical dimension
$\nu(\alpha)$
of a pseudo-effective class
$\alpha \in H^{1,1}(X,\mathbb{R})$
is defined as
\[ \max \{ k \in \{ 0, \ldots, n\}: (\alpha^k \cdot \omega^{n-k})_{\geq 0} 
                                                                     > 0\}\]
for some (and hence all) K\"ahler classes 
$\omega$.
\end{Def}

\noindent A pseudo-effective line bundle $L$ is big iff 
$\nu(L) = \nu(c_1(L)) = n$ (\cite[Thm.3.1.31]{Bou02}). By cutting down with 
ample 
hypersurfaces this shows that the numerical dimension of the first 
Chern class of a pseudo-effective line bundle $L$ is  $\geq \kappa(X,L)$, the 
Kodaira-Iitaka dimension of $L$. 

\noindent Proving that on every projective complex manifold which is not 
uniruled the 
canonical bundle is pseudo-effective (\cite[Cor.0.3]{BDPP04}) Boucksom, 
Demailly, Paun and Peternell were able to use this notion of numerical 
dimension for generalizing the Abundance conjecture to
\begin{conj}
On every projective manifold which is not uniruled we have
\[ \kappa(X) = \nu(X) = \nu(c_1(K_X)). \]
\end{conj}


\noindent The author in turn tried to find more geometric obstacles for 
equality of Kodaira dimension and numerical dimension. In \cite{Eck04b} this 
led him to the notion of numerically trivial foliations. The starting point is
\begin{Def} \label{NumTrivPseff-def}
Let
$X$
be a compact K\"ahler manifold with K\"ahler form 
$\omega$
and pseudo-effective class 
$\alpha \in H^{1,1}(X,\mathbb{R})$.
A submanifold
$Y \subset X$
(closed or not) is \textbf{numerically trivial} w.r.t. 
$\alpha$
iff for every immersed disk 
$\Delta \subset Y$,
\[ \lim_{\epsilon \downarrow 0} \sup_T \int_{\Delta^\prime - \mathrm{Sing}\ T} 
    (T + \epsilon \omega)  = 0, \]
where the
$T$'s
run through all currents with analytic singularities in
$\alpha[-\epsilon \omega]$
and
$\Delta^\prime = \{ t: |t| < 1 - \delta \}$
is any smaller disk contained in 
$\Delta = \{ t: |t| < 1 \}$. 
\end{Def}

\noindent This definition is applied to the leaves of a foliation on $X$ which 
is allowed to have singularities. Such a foliation is given by a saturated 
subsheaf 
\[ \mathcal{F} \subset T_X \]
of the tangent bundle $T_X$ which is closed under the Lie bracket. The 
singularities of $\mathcal{F}$ form the analytic subset $Z$ of points where
\[ \mathcal{F}/m_{X,x}\mathcal{F} \rightarrow T_{X,x} \]
is not injective. By the Frobenius integrability theorem we can cover $X-Z$ by 
open sets
$U_i \cong \Delta^n$
such that there exists smooth holomorphic maps
$p_i: U_i \rightarrow \Delta^{n-k}$
induced by the projection
$\Delta^n \rightarrow \Delta^{n-k}$ with
\[ \mathcal{F}_{|U_i} = T_{U_i/\Delta^{n-k}}. \] 
Further properties of singular foliations and constructions as the union of two
foliations will be discussed in section~\ref{singfol-sec}.
\begin{Def}
Let
$X$
be a compact K\"ahler manifold with a pseudo-effective class
$\alpha \in H^{1,1}(X, \mathbb{R})$. 
A foliation
$\{ \mathcal{F},(U_i,p_i) \}$
is \textbf{numerically trivial} w.r.t.
$\alpha$
iff 
\begin{itemize}
\item[(i)]
every fiber of
$p_i$
is numerically trivial w.r.t. 
$\alpha$,
\item[(ii)]
and if 
$\Delta^2 \hookrightarrow U_i$
is an immersion such that the projection onto the first coordinate coincides 
with the projection
$p_i: U_i \rightarrow \Delta^{n-k}$,
then for any
$\Delta^\prime \subset\subset \Delta$
and any sequence of currents
$T_k \in \alpha[-\epsilon_k\omega]$,
$\epsilon_k \rightarrow 0$, 
the integrals
$\int_{(\{z_1 = a\} \cap \Delta^\prime) - \mathrm{Sing\ }T_k} 
  (T_k + \epsilon_k\omega)$
are uniformly (in
$a$)
bounded from above.
\end{itemize}
\end{Def}

\noindent By proving the Local Key Lemma \cite[Lem.3.8]{Eck04b} the author 
showed that it is possible to construct a numerically trivial foliation maximal
w.r.t. inclusion which is called \textbf{the numerically trivial foliation of 
the pseudo-effective class $\alpha$}. Furthermore, if $\alpha$ is the first 
Chern class $c_1(L)$ of a pseudo-effective line bundle $L$ on $X$ it is shown 
that the fibers of the Kodaira-Iitaka fibration ($m \gg 0$)
\[ \phi_{|mL|}: X \dasharrow Y \]
contain the leaves of the numerically trivial foliation of $L$ 
(\cite[3.3]{Eck04b}). Hence the 
dimension of its leaves is $\leq$ than the dimension of the (generic) fibers of
$\phi_{|mL|}$ and the codimension of the leaves is $\geq$ than
\[ \dim Y = \kappa(X). \] 

\noindent When trying to compare the \textit{numerical} dimension with the 
codimension of the leaves the author discoverded that 
Def.~\ref{NumTrivPseff-def} is not appropriate for this purpose. The point is 
that numerical dimensions are defined via integrals over $n$-dimensional 
complex manifolds whereas Def.~\ref{NumTrivPseff-def} only uses integrals over 
$1$-dimensional disks. Hence the usual difficulties when comparing 
$L^p$-integrable functions for different $p$ occur.
In the end this led the author to change the definition of numerical 
triviality:
\begin{Def} \label{NTFolIntro-Def}
Let $X$ be a compact K\"ahler manifold with K\"ahler form $\omega$ and let
$\alpha \in H^{1,1}(X,\mathbb{R})$ be a pseudo-effective $(1,1)$-class. A 
foliation $\mathcal{F}$ is called numerically trivial w.r.t. $\alpha$ iff for 
all $1 \leq p \leq n-1$ and for all test forms 
$u \in \mathcal{D}^{(n-p,n-p)}(X - \mathrm{Sing}\ \mathcal{F})$ 
 \[ (NT)_u \rule{2cm}{0cm} 
    \lim_{\epsilon \downarrow 0} \sup_{T \in \alpha[-\epsilon\omega]}
       \int_X \left| (T_{ac}+\epsilon\omega)^p \wedge u \right| = 0 
           \rule{5cm}{0cm} \]
where the $T$'s run through all currents with analytic singularities 
representing $\alpha$ with $T \geq -\epsilon\omega$ and $T_{ac}$ is the 
absolute continuous part of $T$ in the Lebesgue decomposition.
\end{Def}

\noindent Here a test form for $\mathcal{F}$ is a smooth $(n-p,n-p)$ form with
compact support outside the singularities of $\mathcal{F}$ whose wedge product
with every $(p,p)$ form in 
$\Lambda^{p,p} \left( T_X^\ast / T_\mathcal{F}^\ast \right)$ is $0$ 
(see section~\ref{NTFol-sec} for details). 

\noindent At least in the surface case this definition of numerical 
triviality is implied by that in \cite{Eck04b}. It is also possible to 
construct a maximal numerical trivial foliation w.r.t. the new definition, see
section~\ref{NTFol-sec}. Furthermore the proof of the follwing theorem becomes 
quite simple:
\begin{thm} \label{AppBd-thm}
Let
$X$
be a compact K\"ahler manifold with K\"ahler form 
$\omega$
and
$\alpha \in H^{1,1}(X,\mathbb{R})$
a pseudo-effective class. Let
$\mathcal{F}$
be the numerically trivial foliation of
$\alpha$. Then the numerical dimension
$\nu(\alpha)$ 
is less or equal to the codimension of the leaves of 
$\mathcal{F}$.
\end{thm}

\noindent And a transversality criterion for detecting numerically trivial 
foliations (see Thm.~\ref{transvers-lem}) allows to show the second aim of 
this note:
\begin{thm} \label{KIfib=NTFol-intro-thm}
Let
$X$
be a K\"ahler manifold and 
$L$
a pseudo-effective line bundle on
$X$.
Suppose that the Kodaira-Iitaka dimension
$\kappa(X,L)$
equals the numerical dimension 
$\nu(X,L)$
of 
$L$.
Then the numerical trivial foliation of
$L$
is the Kodaira-Iitaka fibration of
$L$.
\end{thm}

\noindent 
Both  theorems together imply a sufficient geometric criterion for
\[ \kappa(X,L) < \nu(X,L)  \]
where $L$ is a pseudo-effective line bundle on $X$:  
Suppose that $\mathcal{F}_L$ is the numerically trivial foliation
of $L$ and 
\[ codim \{\mathrm{leaves\ of\ } \mathcal{F}_L\} = \nu(X,L). \]
Then the Kodaira dimension is strictly smaller than the numerical dimension if 
$\mathcal{F}_L$ is a genuine foliation, i.e. 
not induced by a fibration. It would be interesting to know whether the 
converse of this criterion is also true.

\noindent In \cite{Eck04b} 2 surface examples are illustrating this criterion. 
Unfortunately another example \cite[4.3]{Eck04b} dealing with the anticanonical
bundle of $\mathbb{P}^2$ 
blown up in 9 points lying sufficiently general on a smooth elliptic curve 
shows that
\[ codim \{\mathrm{leaves\ of\ } \mathcal{F}_L\} > \nu(X,L) \]
may also occur. Of course, the generalized Abundance Conjecture together with 
Theorem~\ref{KIfib=NTFol-intro-thm} imply that this is not the case for the
canonical bundle.

\section{Singular Foliations} \label{singfol-sec}

\noindent Holomorphic foliations on complex manifolds are usually defined as 
involutive subbundles of the tangent bundle. Then the classical theorem of
Frobenius asserts that through every point there is a unique integral 
complex submanifold~\cite{Miy86}. Singular foliations may be defined as 
involutive coherent subsheaves of the tangent bundle which are furthermore 
saturated, that is their quotient with the tangent bundle is torsion free.
In points where the rank is maximal one may use again the Frobenius 
theorem to get leaves.

\noindent Later on we use the following notation:
\begin{Def}
Let $X$ be a complex manifold and $\mathcal{F} \subset T_X$
a saturated involutive subsheaf. Then the analytic subset
\[ \left\{ x \in X : 
   \mathcal{F}/m_{X,x}\mathcal{F}  \rightarrow T_{X,x} \ 
   \mathrm{is\ not\ injective} \right\}\]
is called the singular locus of $\mathcal{F}$ and is denoted by 
$\mathrm{Sing\ } \mathcal{F}$. The dimension of 
$\mathcal{F}/m_{X,x}\mathcal{F}$ in a point 
$x \in X - \mathrm{Sing\ } \mathcal{F}$ is called rank of $\mathcal{F}$ and
denoted by $\mathrm{rk}(\mathcal{F})$.
\end{Def}

\noindent Because $\mathcal{F}$ is saturated we have
$\mathrm{codim}\ \mathrm{Sing\ } \mathcal{F} \geq 2$.
The existence of leaves means that around every point 
$x \in X - \mathrm{Sing\ } \mathcal{F}$ there is an (analytically) open subset 
$U \subset X - \mathrm{Sing\ } \mathcal{F}$ with
coordinates $z_1, \ldots z_n$, $n = \dim X$, such that the leaves of 
$\mathcal{F}$ are the fibers of the projection onto the coordinates
$z_{k+1}, \ldots, z_n$ where $k = \mathrm{rk}(\mathcal{F})$. In particular the 
leaves have dimension $\mathrm{rk}(\mathcal{F})$.

\noindent To construct numerically trivial foliations we need a local 
description of several 
operations applied on two foliations. We start with the easiest configuration:
\begin{prop} \label{GsubF-prop}
Let $\mathcal{G} \subset \mathcal{F}$ be two foliations on a complex manifold 
$X$, $\mathrm{rk}(\mathcal{F}) = k$, $\mathrm{rk}(\mathcal{G}) = l$, $l < k$.
Then for all $x \in X - (\mathrm{Sing\ } \mathcal{F} \cup 
                         \mathrm{Sing\ } \mathcal{G})$ 
there is an open neighborhood $U \subset X - (\mathrm{Sing\ } \mathcal{F} \cup 
                                              \mathrm{Sing\ } \mathcal{G})$
with coordinates $z_1, \ldots, z_n$ such that the leaves of $\mathcal{F}$ are
the fibers of the projection onto the last $n-k$ coordinates and the leaves of 
$\mathcal{G}$ are
the fibers of the projection onto the last $n-l$ coordinates.
\end{prop}
\begin{proof}
This is an easy consequence of the theorem on implicitely defined functions. 
Note that neither $\mathrm{Sing\ } \mathcal{G}$ need to be contained in 
$\mathrm{Sing\ } \mathcal{F}$ nor vice versa.
\end{proof}

\begin{Def}
Let $\mathcal{F}$ and $\mathcal{G}$ be two foliations on a complex manifold 
$X$. Then $\mathcal{F} \cap \mathcal{G} \subset T_X$ is called the intersection
foliation of 
$\mathcal{F}$ and $\mathcal{G}$.
\end{Def}

\noindent Note that $\mathcal{F} \cap \mathcal{G}$ is certainly involutive but 
may be not saturated: the rank of $\mathcal{F} \cap \mathcal{G}$ can even jump
in $\mathrm{codim}\ 1$ subsets. To get a better picture in local coordinates we
nevertheless think of it as a foliation and denote by 
$\mathrm{Sing\ } (\mathcal{F} \cap \mathcal{G})$ the analytic locus where the
rank jumps.

\begin{prop} \label{FGcoord-prop}
Let $\mathcal{F}$ and $\mathcal{G}$ be two foliations on a complex manifold $X$
with $\mathrm{rk}(\mathcal{F}) = k$, $\mathrm{rk}(\mathcal{G}) = m$ and 
$\mathrm{rk}(\mathcal{F} \cap \mathcal{G}) = l$. Let $x \in X$ be a point which
is not singular for $\mathcal{F}$, $\mathcal{G}$ and 
$\mathcal{F} \cap \mathcal{G}$. Then there exists an open neighborhood
\[ U \subset X - (\mathrm{Sing\ } \mathcal{F} \cup \mathrm{Sing\ } \mathcal{G}
\cup \mathrm{Sing} \left(\mathcal{F} \cap \mathcal{G}\right)) \]
of $x$ with coordinates $z_1, \ldots, z_n$ such that 
\begin{itemize}
\item[(i)] the leaves of $\mathcal{F}$ in $U$ are the fibers of the projection 
on $z_{k+1}, \ldots , z_n$,
\item[(ii)] the leaves of $\mathcal{F} \cap \mathcal{G}$ in $U$ are the fibers 
  of the projection on $z_{l+1}, \ldots , z_n$ and
\item[(iii)] the leaves of $\mathcal{G}$ in $U$ are the fibers of the 
  projection on $z_{l+1}, \ldots , z_k, g_{m+k-l+1}(z), \ldots, g_n(z)$
  where the $g$'s are analytic functions with
  $g_{m+k-l+j}(z)_{|U_x} = z_{k+j}$ on 
  \[ U_x = \{ z \in U : z_{l+1}(z) = z_{l+1}(x), \ldots, z_k(z) = z_k(x) \}. \]
\end{itemize}
\end{prop} 
\begin{proof}
Again this results from applying the theorem on implicitely defined functions
several times. Since the geometry is more difficult than in 
Prop.~\ref{GsubF-prop} (see Figure~\ref{UnionFol-fig}) we present more details:
Choose coordinates $z_1, \ldots, z_n$ for $\mathcal{F}$ and 
$\mathcal{F} \cap \mathcal{G}$ in a neighborhood
\[ U^\prime \subset X - (\mathrm{Sing\ } \mathcal{F} \cup 
                         \mathrm{Sing\ } \mathcal{G} \cup 
                         \mathrm{Sing\ } (\mathcal{F} \cap \mathcal{G})) \]
of $x$ as in Prop.~\ref{GsubF-prop}. Since the leaves of $\mathcal{G}$ contain
the leaves of $\mathcal{F} \cap \mathcal{G}$ we can describe the leaves of
$\mathcal{G}$ in $U^\prime$ (possibly restricted) as the fibers of the 
projection given by analytic functions $g_{m+1}, \ldots, g_n$ only depending on
$z_{l+1}, \ldots , z_n$. Furthermore we know that for a fixed point 
$(z_{k+1}, \ldots, z_n)$ the fibers of the projection from the leaf of
$\mathcal{F}$ given by $g_{m+1}, \ldots, g_n$ are leaves of 
$\mathcal{F} \cap \mathcal{G}$.

\noindent Consequently an applicaton of the theorem on implicitely defined 
functions gives us (after possibly reordering the $g$'s) coordinates
\[ z_1, \ldots, z_l, z_{l+1}^\prime = g_{m+1}(z), \ldots, 
                     z_{k}^\prime = g_{m+k-l}(z), z_{k+1}, \ldots, z_n \]
in an open subset $U^{\prime\prime} \subset V^\prime$
such that the leaves of $\mathcal{F}$ resp. $\mathcal{F} \cap \mathcal{G}$ are
still the projection onto the last $n-k$ resp. $n-l$ coordinates, and the 
leaves of $\mathcal{G}$ are the fibers of the projection onto
\[ z_{l+1}^\prime, \ldots, z_{k}^\prime, g_{m+k-l+1}(z), \ldots , g_n(z). \]

\noindent Now we fix $z_{l+1}^\prime = a_{l+1}, \ldots, z_k^\prime = a_k$.
Using again the theorem on implicitely defined functions we see that
(after possibly another reordering of the $g$'s) 
\[ \begin{array}{l}
   z_1, \ldots, z_l, z_{l+1}^\prime , \ldots, z_{k}^\prime , \\
   z_{k+1}^\prime = g_{m+k-l+1}(a_{l+1}, \ldots, a_k, z_{k+1}, \ldots, z_n),\\ 
   \vdots \\ 
   z_{n-m+l}^\prime = g_{n}(a_{l+1}, \ldots, a_k, z_{k+1}, \ldots, z_n), \\
   z_{n-m+l+1}, \ldots, z_n 
   \end{array} \]
are coordinates in an open subset $U \subset U^{\prime\prime}$ having all the
properties claimed in the proposition.
\end{proof}

\noindent For our purposes the most important operation on two holomorphic 
foliations $\mathcal{F}$ and $\mathcal{G}$ on a complex manifold $X$ is the 
\textit{union} $\mathcal{F} \sqcup \mathcal{G}$. We define it as the foliation
given by the smallest 
saturated involutive subsheaf of $T_X$ containing both $\mathcal{F}$ and 
$\mathcal{G}$. Such a sheaf exists because saturated foliations contained in 
each other have different ranks, the intersection of two foliations is again a 
foliation and $T_X$ is involutive.

\noindent Besides this pure existence statement there is an inductive algebraic
construction of $\mathcal{F} \sqcup \mathcal{G}$:
\[ \begin{array}{rcl}
   \mathcal{H}_1 & := & \mathrm{saturation\ of\ } \mathcal{F} + \mathcal{G} \\
   \mathcal{H}_2 & := & \mathrm{saturation\ of\ } \mathcal{H}_1 + 
                                          [\mathcal{H}_1,\mathcal{H}_1] \\
    & \vdots &
   \end{array} \]
and so on until $\mathcal{H}_m = \mathcal{H}_{m+1}$ which means 
$[\mathcal{H}_m,\mathcal{H}_m] \subset \mathcal{H}_m$. Then 
$\mathcal{H}_m = \mathcal{F} \sqcup \mathcal{G}$. This is a local construction 
hence for open subsets $U \subset X$ we have
\[ \mathcal{F}_{|U} \sqcup \mathcal{G}_{|U} = 
   (\mathcal{F} \sqcup \mathcal{G})_{|U}. \]

\noindent We want to describe an inductive geometric construction of 
$\mathcal{F} \sqcup \mathcal{G}$ on open subsets
\[ U \subset X - (\mathrm{Sing\ } \mathcal{F} \cup \mathrm{Sing\ } \mathcal{G}
                  \cup \mathrm{Sing\ } (\mathcal{F} \cap \mathcal{G})) - Z \]
where $Z$ is an analytic subset of 
$X - (\mathrm{Sing\ } \mathcal{F} \cup \mathrm{Sing\ } \mathcal{G}
                  \cup \mathrm{Sing\ } (\mathcal{F} \cap \mathcal{G}))$ which 
will be determined during the construction. 
Following the inductive steps of this construction we will later on prove
Key Lemma~\ref{KeyLemma}.

\noindent Start with a neighborhood $U$ of a point 
$x \in X - (\mathrm{Sing\ } \mathcal{F} \cup \mathrm{Sing\ } \mathcal{G}
                  \cup \mathrm{Sing\ } (\mathcal{F} \cap \mathcal{G}))$
having coordinates $z_1, \ldots, z_n$ as in Prop.~\ref{FGcoord-prop}. Define a
foliation $\mathcal{G}^\prime$ on $U$ whose leaves are the fibers of the 
projection on $z_{l+1}, \ldots , z_{n-m+l}$. Figure~\ref{UnionFol-fig} 
illustrates that in general
$\mathcal{F} + \mathcal{G}^\prime \neq \mathcal{F} \sqcup \mathcal{G}$ (take 
the fibers of the vertical projection as leaves of $\mathcal{F}$ whereas the 
leaves of $\mathcal{G}$ are the horizontal lines twisted around in vertical 
direction):

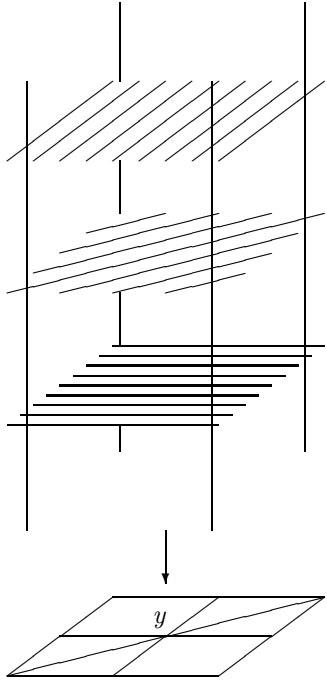
\begin{figure} \label{UnionFol-fig}
\begin{center}
\begin{picture}(130,260)(0,0)
\put(7.5,55){\line(0,1){170}}
\put(77.5,55){\line(0,1){170}}
\put(42.5,85){\line(0,1){10}}
\put(42.5,125){\line(0,1){20}}
\put(42.5,175){\line(0,1){20}}
\put(42.5,225){\line(0,1){30}}
\put(112.5,85){\line(0,1){170}}

\put(0,95){\line(1,0){80}}
\put(5,98.7){\line(1,0){80}}
\put(10,102.5){\line(1,0){80}}
\put(15,106.2){\line(1,0){80}}
\put(20,110){\line(1,0){80}}
\put(25,113.7){\line(1,0){80}}
\put(30,117.5){\line(1,0){80}}
\put(35,121.2){\line(1,0){80}}
\put(40,125){\line(1,0){80}}

\put(0,195){\line(4,3){40}}
\put(10,195){\line(4,3){40}}
\put(20,195){\line(4,3){40}}
\put(30,195){\line(4,3){40}}
\put(40,195){\line(4,3){40}}
\put(50,195){\line(4,3){40}}
\put(60,195){\line(4,3){40}}
\put(70,195){\line(4,3){40}}
\put(80,195){\line(4,3){40}}

\put(0,145){\line(4,1){120}}
\put(20,145){\line(4,1){90}}
\put(40,145){\line(4,1){60}}
\put(60,145){\line(4,1){30}}
\put(100,175){\line(-4,-1){90}}
\put(80,175){\line(-4,-1){60}}
\put(60,175){\line(-4,-1){30}}

\put(60,55){\vector(0,-1){20}}

\put(0,0){\line(1,0){80}}
\put(40,30){\line(1,0){80}}
\put(0,0){\line(4,3){40}}
\put(80,0){\line(4,3){40}}

\put(55,20){$\mathit{y}$}

\put(20,15){\line(1,0){80}}
\put(0,0){\line(4,1){120}}
\put(40,0){\line(4,3){40}}
 
\end{picture}
\end{center}
\caption{Union of two foliations}

\end{figure}

\noindent Denoting the projection on $z_{k+1}, \ldots, z_n$ by 
$\pi_\mathcal{F}$ we examine instead $r$-tuples of points $x_1, \ldots, x_r$ in
fibers $\pi_\mathcal{F}^{-1}(y)$ of points 
$y \in \pi_\mathcal{F}(U) \subset \mathbb{C}^{n-k}$. If 
$T_\mathcal{G}(x_i) \subset T_{X,x_i}$ indicates the space of directions 
tangent to $\mathcal{G}$ in $x_i$ we have a sequence of inclusions
\[ 0 \subset d\pi_\mathcal{F}(T_\mathcal{G}(x_1)) \subset 
   d\pi_\mathcal{F}(T_\mathcal{G}(x_1)) + d\pi_\mathcal{F}(T_\mathcal{G}(x_2))
   \subset \cdots \subset \sum_{i=1}^r d\pi_\mathcal{F}(T_\mathcal{G}(x_i))
   \subset T_{\mathbb{C}^{n-k},y}. \]
There is an $r \in \mathbb{N}$ and a Zariski open subset of the $r$-fold 
product
\[ \pi_\mathcal{F}^{-1}(y) \times \cdots \times \pi_\mathcal{F}^{-1}(y) \]
such that
\begin{itemize}
\item[(i)] all inclusions in the above sequence are strict and
\item[(ii)] $d\pi_\mathcal{F}(T_\mathcal{G}(x^\prime)) \subset 
             \sum_{i=1}^r d\pi_\mathcal{F}(T_\mathcal{G}(x_i))$ for every point
$x^\prime \in \pi_\mathcal{F}^{-1}(y)$.
\end{itemize}

\noindent Varying $y \in \pi_\mathcal{F}(U)$ may change the number $r$ and the 
dimensions of the vector spaces
\[ \sum_{i=1}^s d\pi_\mathcal{F}(T_\mathcal{G}(x_i)),\ s = 1, \ldots, r. \]
But again there is an analytic subset $Z_U \subset \pi_\mathcal{F}(U)$ such 
that for $y \in V := \pi_\mathcal{F}(U) - Z_U$ the dimensions and $r$ remain 
constant. Since everything is defined intrinisically the sets 
$\pi_\mathcal{F}^{-1}(Z_U)$ glue together to an analytic subset $Z$ of 
$X - (\mathrm{Sing\ } \mathcal{F} \cup \mathrm{Sing\ } \mathcal{G}
                  \cup \mathrm{Sing\ } (\mathcal{F} \cap \mathcal{G}))$.
Furthermore we can find $r$ sections $\sigma_i : V^\prime \rightarrow U$ of
$\pi_\mathcal{F}$ such that
\begin{itemize}
\item[(i)] the points $x_i := \sigma_i(y)$ produce a sequence of tangent 
subspaces 
$\sum_{i=1}^s d\pi_\mathcal{F}(T_\mathcal{G}(x_i))$, $s = 1, \ldots, r$, 
strictly included in each other and
\item[(ii)] if $\pi: U \rightarrow \mathbb{C}^{k-l}$ is the projection onto 
$z_{l+1}, \ldots, z_k$, the map $\pi \circ \sigma_i$ is constant.
\end{itemize}
The $V^\prime$ are possibly smaller subsets of $V$ covering $V$.

\noindent To get the announced inductive construction of 
$\mathcal{F} \sqcup \mathcal{G}$ on $\pi_\mathcal{F}^{-1}(V^\prime)$ we need 
another little observation: Since the holomorphic functions $g_j$ defining 
$\pi_\mathcal{G}$ do not depend on $z_1, \ldots, z_l$ (see proof of 
Prop.~\ref{FGcoord-prop}) the tangent space
\[ d\pi_\mathcal{F}(T_\mathcal{G}(x)) \]
does not change for different $x$ in the intersection of a fixed 
$\pi_\mathcal{F}$- and a $\pi$-fiber. Furthermore the fibers of $\pi$ consist 
of leaves of $\mathcal{G}$.

\noindent Now we construct inductively foliations $\mathcal{F}_i$, 
$i = 0, \ldots, r$ on $\pi_\mathcal{F}^{-1}(V^\prime)$. We start with
\[ \mathcal{F}_0 := \mathcal{F} \cap \pi_\mathcal{F}^{-1}(V^\prime). \]
Because of the observation above the leaves of $\mathcal{G}$ in 
$\pi^{-1}(\pi(x_1))$ map onto the leaves of a smooth foliation $\mathcal{G}_1$
on $V^\prime$ which is induced by a projection $\pi_{\mathcal{G}_1}$. Put
\[ \mathcal{F}_1 := \pi_\mathcal{F}^{-1}(\mathcal{G}_1) \]
and let $\pi_{\mathcal{F}_1} := \pi_{\mathcal{G}_1} \circ \pi_\mathcal{F}$ be 
the projection whose fibers are the leaves of $\mathcal{F}_1$.

\noindent The observation and the properties of the $x_1, \ldots, x_r$ imply 
that $T_{\mathcal{G}|\pi^{-1}(\pi(x_2))}$ maps onto an involutive 
\textbf{subbundle} of $T_{\pi_{\mathcal{F}_1}(\pi_\mathcal{F}^{-1}(V^\prime))}$
and consequently the leaves of $\mathcal{G}$ in $\pi^{-1}(\pi(x_2))$ also map
onto leaves of a smooth foliation $\mathcal{G}_2$ on 
$\pi_{\mathcal{F}_1}(\pi_\mathcal{F}^{-1}(V^\prime))$. Define
\[ \mathcal{F}_2 := \pi_{\mathcal{F}_1}^{-1}(\mathcal{G}_2) \]
and continue inductively setting
\[ \mathcal{F}_i := \pi_{\mathcal{F}_{i-1}}^{-1}(\mathcal{G}_i) \]
where $\mathcal{G}_i$ is the image of the leaves of $\mathcal{G}$ in 
$\pi^{-1}(\pi(x_i)$ on 
$\pi_{\mathcal{F}_{i-1}}(\pi_\mathcal{F}^{-1}(V^\prime))$.

\noindent By construction these foliations $\mathcal{F}_s$ have as tangent 
space in a point $x \in \pi_\mathcal{F}^{-1}(V^\prime)$
\[ d\pi_\mathcal{F}(x)^{-1}(\sum_{i=1}^s d\pi_\mathcal{F}(T_\mathcal{G}(x_i))\]
where $\pi_\mathcal{F}(x_i) = \pi_\mathcal{F}(x)$ for all $i$.
In addition $\mathcal{F}_r$ contains all leaves of $\mathcal{F}$ and 
$\mathcal{G}$ in $\pi_\mathcal{F}^{-1}(V^\prime)$: otherwise there is a point 
$y \in V^\prime$ and a point $x \in \pi_\mathcal{F}^{-1}(y)$ such that
\[ d\pi_\mathcal{F}(T_\mathcal{G}(x)) \not\subset
   \sum_{i=1}^r d\pi_\mathcal{F}(T_\mathcal{G}(x_i)), \]
$\pi_\mathcal{F}(x_i) = \pi_\mathcal{F}(x)$ for all $i$.

\noindent On the other hand $T_{\mathcal{F} \sqcup \mathcal{G}}(x)$ must 
contain every tangent subspace
\[ d\pi_\mathcal{F}(x)^{-1}(T_\mathcal{G}(x^\prime)) \]
of points $x^\prime$ with $\pi_\mathcal{F}(x^\prime) = \pi_\mathcal{F}(x)$
since 
$\pi_\mathcal{F}^{-1}\pi_\mathcal{F}
 (\pi_\mathcal{G}^{-1}\pi_\mathcal{G}(x^\prime))$ is contained in a leaf of 
$\mathcal{F} \sqcup \mathcal{G}$. Consequently, 
$d\pi_\mathcal{F}(x)^{-1}(\sum_{i=1}^r d\pi_\mathcal{F}(T_\mathcal{G}(x_i))
 \subset T_{\mathcal{F} \sqcup \mathcal{G}}(x)$ and on 
$\pi_\mathcal{F}^{-1}(V^\prime)$ we have
\[ \mathcal{F} \sqcup \mathcal{G} = \mathcal{F}_r . \]

\noindent An important type of foliations are those induced in a unique way by
meromorphic maps $f: X \dasharrow Y$ from a compact complex manifold $X$ to 
another complex manifold $Y$: Take the relative tangent sheaf of $f$ on the 
Zariski open subset $U$ where $f$ is smooth and saturate.
\begin{prop} \label{fibjoin-prop}
Let $X$ be a compact complex manifold and $f: X \dasharrow Y_1$, 
$g: X \dasharrow Y_2$ two meromorphic maps with induced foliations
$\mathcal{F}$ and $\mathcal{G}$ on $X$. Then $\mathcal{F} \sqcup \mathcal{G}$ 
is also induced by a meromorphic map $h: X \dasharrow Z$. 
\end{prop}
\begin{proof}
Let $\Gamma_f \subset Y_1 \times X$, $\Gamma_g \subset X \times Y_2$ be the 
graphs of $f$ and $g$. Consider the product $Y_1 \times X \times X \times Y_2$
and its projections $p_1, p_2, p_3, p_4$ onto the subsequent factors. A general
point $(y_1,x_1,x_2,y_2)$ of the intersection
\[ (p_1 \times p_2)^{-1}(\Gamma_f) \cap 
   (p_2 \times p_3 \times p_4)^{-1}(\Gamma_g \times_{Y_2} \Gamma_g) \subset 
   Y_1 \times X \times X \times Y_2 \]
satisfies $x_1 \in f^{-1}(y_1)$, $x_1,x_2 \in g^{-1}(y_2)$. Since $y_1$ and 
$y_2$ are uniquely determined by a general point $x_1 \in X$ there is a unique
irreducible component $W \subset Y_1 \times X \times X \times Y_2$ in this 
intersection projecting surjectively on $X$ via $p_2$ such that the fiber over 
a general point $x \in X$ is a unique $g$-fiber. 
$(p_1 \times p_3)(W) \subset Y_1 \times X$ is a family of compact complex 
subspaces of $X$ parametrized by $Y_1$ and covering $X$. Hence there is a 
meromorphic map from $Y_1$ to the Douady space of $X$ and we call the image of
this map $Z_1$. The image $W_1 \subset Z_1 \times X$ of $(p_1 \times p_3)(W)$ 
in the universal family over the Douady space of $X$ has the following 
properties:
\begin{itemize}
\item[(i)] Every fiber of $W_1$ over $Z_1$ consists of $g$-fibers.
\item[(ii)] For two points in the same fiber of $W_1$ over $Z_1$ there is a 
sequence of $f$- and $g$-fibers connecting them.
\item[(iii)] Given an analytic subset $Z \subset X$ two general points in a 
general fiber may be connected by a 
sequence of $f$- and $g$-fibers such that two subsequent fibers do not 
intersect in $Z$.
\end{itemize}
Points connected by a sequence of $f$- and $g$-fibers satisfying (iii) with
\[ Z = \mathrm{Sing}(\mathcal{F} \sqcup \mathcal{G}) \cup
       \{\mathrm{indeterminacy\ loci\ of\ }f\ \mathrm{and\ }g\} \]
must lie in the same leaf of $\mathcal{F} \sqcup \mathcal{G}$.

\noindent Repeat the construction from above but interchange the r\^oles of $f$
and $g$: Take the generic irreducible component of the intersection
\[ (p_1 \times p_2 \times p_3)^{-1}(\Gamma_f \times_{Y_1} \Gamma_f) \cap
   (p_3 \times p_4)^{-1}(\Gamma_g) \subset Y_1 \times X \times X \times Y_2 \]
under the projection $p_3$, project it via $p_2 \times p_4$ in $X \times Z_1$
and map it via the universal properties of the Douady space into a family 
$W_2 \subset Z_2 \times X$ of compact complex subspaces covering $X$. Now every
fiber of $W_2$ consists of $f$-fibers and $W_2$ also satisfies (ii) and (iii).

\noindent Continue this construction interchanging the r\^oles of $f$ and $g$
in each step until the fiber dimension of $W_k$ over $Z_k$ does not rise in the
next step. This is the case iff for every fiber $F$ and every point $x \in F$
there is an $f$- and $g$-fiber through $x$ contained in $F$. Through a general 
point $x$ these fibers are unique. Since $W_k$ also satisfies (iii) two fibers 
containing the same general point $x$ must be equal. Hence 
$W_k \subset Z_k \times X \rightarrow X$ is generically 1:1 and defines a
meromorphic map $h: X \dasharrow Z_k$ whose induced foliation $\mathcal{H}$
contains 
$\mathcal{F}$ and $\mathcal{G}$. On the other hand the fibers of $W_k$ are 
contained in leaves of $\mathcal{F} \sqcup \mathcal{G}$ by construction. Taken 
all together we have shown: 
\[ \mathcal{H} = \mathcal{F} \sqcup \mathcal{G}. \]
\end{proof}

\begin{rem}
This construction closely resembles that of Campana's reduction map 
\cite{Cam81, Cam94}. The difference is that in the above construction for any 
given analytic subset $Z \subset X$ two general points lie in the same fiber 
iff they can be connected without touching $Z$. The easiest example where the 
two reduction maps fall apart are the quotients with respect to the pencil of
lines through a point in $\mathbb{P}^2$.
\end{rem}

\noindent We generalize this construction to the following
\begin{Def}
Let $X$ be a compact complex manifold. A covering family 
$(C_t)_{t \in T}$ of complex subspaces of $X$ is called generically 
connecting iff for any analytic subset $Z \subset X$ two general points are 
connected by a finite sequence of elements in $(C_t)$ such that two subsequent 
elements do not intersect in $Z$.

\noindent A meromorphic map $f: X \dasharrow Y$ is called the generic reduction
map with respect to $(C_t)_{t \in T}$ iff the general fibers are generically
$C_t$-connected and every element of $(C_t)$ is contained in a fiber. Here, 
fibers of $f$ are defined via the graph of $f$.
\end{Def} 

\noindent The construction described above shows that for every family 
$(C_t)_{t \in T}$ there exists a unique generic reduction map.

\noindent Another difference between the generic and Campana's reduction map is
the stability under modifications: Let $X$ be a compact complex manifold and 
$(C_t)_{t \in T}$ a covering family of complex subspaces of $X$. Let 
$f : X \dasharrow Y$ be the generic quotient and $g : X \dasharrow Z$ Campana's
quotient with respect to $(C_t)$. If $\pi:\widehat{X} \rightarrow X$ is a 
modification of compact K\"ahler manifolds then the generic quotient of
$\widehat{X}$ w.r.t. the strict or total transforms of $(C_t)$ is 
described by $f \circ \pi$ whereas in general Campana's quotient is described
by $g \circ \pi$ only w.r.t. the total transforms of $(C_t)$ (cf. the pencil of
lines through a point in $\mathbb{P}^2$).

\section{Numerically Trivial Foliations} \label{NTFol-sec}

\noindent From now on let $X$ be a compact K\"ahler manifold. On 
$X - \mathrm{Sing}\ \mathcal{F}$ a foliation $\mathcal{F}$ is described by a 
subbundle $T_\mathcal{F}$ of the tangent bundle $T_X$. Then 
$(T_X/T_\mathcal{F})^\ast$ is a (holomorphic) subbundle of 
$T_X^\ast = \Omega_X^1$. Hence the subbundle generates $(p,p)$-forms on 
$X - \mathrm{Sing}\ \mathcal{F}$ which we collect in the set
\[ \mathcal{E}^{p,p}(X, \mathcal{F}) \subset 
   \mathcal{E}^{p,p}(X - \mathrm{Sing}\ \mathcal{F}). \]
Note that $\mathcal{E}^{p,p}(X, \mathcal{F}) = 0$ for 
$p > \dim X - \mathrm{rk\ } \mathcal{F}$.

\noindent An $(n-p,n-p)$-form 
$u \in \mathcal{D}^{(n-p,n-p)}(X - \mathrm{Sing}\ \mathcal{F})$ is called 
\textbf{test form} for $\mathcal{F}$ iff for all 
$v \in \mathcal{E}^{p,p}(X, \mathcal{F})$
\[ v \wedge u = 0. \]
\begin{Def} \label{NTFol-Def}
Let $X$ be a compact K\"ahler manifold with K\"ahler form $\omega$ and let
$\alpha \in H^{1,1}(X,\mathbb{R})$ be a pseudo-effective $(1,1)$-class. A 
foliation $\mathcal{F}$ is called numerically trivial w.r.t. $\alpha$ iff for 
all $1 \leq p \leq n-1$ and for all test forms 
$u \in \mathcal{D}^{(n-p,n-p)}(X - \mathrm{Sing}\ \mathcal{F})$ 
 \[ (NT)_u \rule{2cm}{0cm} 
    \lim_{\epsilon \downarrow 0} \sup_{T \in \alpha[-\epsilon\omega]}
       \int_X \left| (T_{ac}+\epsilon\omega)^p \wedge u \right| = 0 
           \rule{5cm}{0cm} \]
where the $T$'s run through all currents with analytic singularities 
representing $\alpha$ with $T \geq -\epsilon\omega$ and $T_{ac}$ is the 
absolute continuous part of $T$ in the Lebesgue decomposition.
\end{Def}

\noindent Note that $\left| (T_{ac}+\epsilon\omega)^p \wedge u \right|$ is the
total variation of the measure $(T_{ac}+\epsilon\omega)^p \wedge u$. Since all
occuring currents are absolutely continuous they may be written as forms with 
(at least) $L^1_{loc}$ functions as coefficients. In particular 
$(T_{ac}+\epsilon\omega)^p \wedge u = f \cdot \omega^n$ for an 
$L^1_{loc}$ function $f$ and
\[ \left| (T_{ac}+\epsilon\omega)^p \wedge u \right| = |f| \cdot \omega^n. \]
The usual facts about absolute values like the triangle inequality follow 
immediately from this formula.

\noindent To verify numerical triviality of a foliation $\mathcal{F}$ we only 
need to check condition $(NT)_u$ for special test forms for $\mathcal{F}$:
\begin{Def}
Let $\mathcal{F}$ be a foliation of rank $k$ on a complex $n$-dimensional
manifold $X$ and $U \subset X - \mathrm{Sing}\ \mathcal{F}$ an open subset with
coordinates $z_1, \ldots, z_n$ such that the leaves of $\mathcal{F}$ in $U$ are
the fibers of the projection onto the last $n-k$ coordinates.

\noindent Let $1 \leq p \leq n-1$. An $(n-p,n-p)$ form 
$\sum_{|I|=|J|=n-p} a_{IJ} dz_I \wedge d\overline{z}_J$ on $U$ is called  
constant test form for $\mathcal{F}$ on $U$ iff the $a_{IJ} \in \mathbb{C}$
are constant and
\[ \#(I \cap \{k+1,\ldots,n\}) \leq n-k-p\ \mathrm{and\ }
   \#(J \cap \{k+1,\ldots,n\}) \leq n-k-p\ \ \Longrightarrow a_{IJ} = 0. \]  
\end{Def}

\begin{thm} \label{NTCrit-thm}
Let $X$ be a compact K\"ahler manifold with K\"ahler form $\omega$ and let
$\alpha \in H^{1,1}(X, \mathbb{R})$ be a pseudoeffective $(1,1)$ class. A 
foliation $\mathcal{F}$ on $X$ of $\mathrm{rk\ }k$ is numerically trivial 
w.r.t. $\alpha$ iff there exists a covering $\{ U_i \}$ of relatively compact
open subsets of $X - \mathrm{Sing\ }\mathcal{F}$ satisfying the following 
conditions: On each $U_i$ there are coordinates $z_1, \ldots, z_n$ such that
\begin{itemize}
\item[(i)]
the leaves of $\mathcal{F}$ are the fibers of the projection onto the last 
$n-k$ coordinates, 
\item[(ii)] for all $1 \leq p \leq n-1$ and for all real constant test forms
$u$ for $\mathcal{F}$ in $U_i$ the equality $(NT)_u$ is true.
\end{itemize}
\end{thm}
\begin{proof}
In each point $x \in U_i$ constant test forms for $\mathcal{F}$ in $U_i$ 
generate all forms of $\mathcal{E}^{p,p}(X, \mathcal{F})$. Hence arbitrary test
forms for $\mathcal{F}$ in $U_i$ can be approximated by locally constant test
forms: Let $u$ be a real $(n-p,n-p)$ test form for $\mathcal{F}$ in $U_i$. Then
for every $\epsilon > 0$ there exists a locally constant test form $u_\epsilon$
such that
\[ -\epsilon\omega^{n-p} < u - u_\epsilon < \epsilon\omega^{n-p} \]  
in $U_i$. Consequently we get 
\begin{eqnarray*} 
    \int_{U_i} \left| (T_{ac}+\epsilon\omega)^p \wedge u \right| & \leq &  
    \int_{U_i} \left| (T_{ac}+\epsilon\omega)^p \wedge (u-u_\epsilon ) \right|
  + \int_{U_i} \left| (T_{ac}+\epsilon\omega)^p \wedge u_\epsilon \right| \\
     & \leq & \epsilon \int_{U_i} (T_{ac}+\epsilon\omega)^p \wedge \omega^{n-p}
  + \int_{U_i} \left| (T_{ac}+\epsilon\omega)^p \wedge u_\epsilon \right| 
\end{eqnarray*} 
for every current $T \in \alpha[-\epsilon\omega]$ with analytic singularities.
By Boucksom's theory of moving intersection numbers \cite{Bou02}[3.2] we have
\[ \lim_{\epsilon \downarrow 0} \sup_T 
   \int_{U_i} (T_{ac}+\epsilon\omega)^p \wedge \omega^{n-p} \leq 
   \lim_{\epsilon \downarrow 0} \sup_T 
   \int_X (T_{ac}+\epsilon\omega)^p \wedge \omega^{n-p} = 
   (\alpha.\omega^{n-p})_{\geq 0}. \]

\noindent For the second summand note that $u_\epsilon$ is a linear combination
of globally constant test forms for $\mathcal{F}$ in $U_i$ multiplied with 
characteristic functions. Using the assumption we conclude
\[ \lim_{\epsilon \downarrow 0} \sup_T 
   \int_{U_i} \left| (T_{ac}+\epsilon\omega)^p \wedge u_\epsilon \right| = 0.\]
\end{proof}

\noindent To justify the definition of numerical triviality we prove
\begin{lem} \label{pullback-lem}
Let $U \subset \mathbb{C}^n$ be an open subset and 
$\pi: U \rightarrow \mathbb{C}^{n-k}$ the projection onto the last $n-k$ 
factors. If $T$ is a closed absolutely continuous $(1,1)$-current on $U$ such
that
\[ \int_U \left| T \wedge u \right| = 0 \]
for all test forms $u \in \mathcal{D}^{n-1,n-1}(U)$ for $\pi$ then there exists
a $(1,1)$-current $S$ on $\pi(U) \subset \mathbb{C}^{n-k}$ such that
\[ T \equiv \pi^\ast S. \]
\end{lem}
\begin{proof}
Since $T$ is absolutely continuous there are $L^1_\mathrm{loc}$-functions 
$f_{ij}$ such that
\[ T \equiv \sum f_{ij} dz_i \wedge d\overline{z}_j. \]
Hence $\int_U \left| T \wedge u \right| = 0$ for all $\pi$-test forms
$u \in \mathcal{D}^{n-1,n-1}(U)$ tells us  
\[ (i,j) \not\in \{ k+1, \ldots, n \} \times \{ k+1, \ldots, n \} 
   \Longrightarrow f_{ij} \equiv 0. \] 
But then the closedness of $T$ implies that for other $f_{ij}$'s the partial 
derivatives $\frac{\partial}{\partial z_l}f_{ij}$ and 
$\frac{\partial}{\partial \overline{z}_l}f_{ij}$ vanish (in the sense of 
currents) if
$l \not\in \{ k+1, \ldots, n \}$. Consequently these $f_{ij}$'s do not depend 
on $z_l$, $l = 1, \ldots, k$, and 
\[ T \equiv \sum_{i,j = k+1}^n f_{ij} dz_i \wedge d\overline{z}_j \]
may be interpreted as pulled back from a current $S$ on $\pi(U)$.
\end{proof}

\noindent We want to show that there is always a \textit{maximal} numerically
trivial foliation $\mathcal{F}$ for a pseudoeffective $(1,1)$-class $\alpha$ 
meaning that every foliation numerically trivial w.r.t. $\alpha$ is contained
in $\mathcal{F}$. This will be a direct consequence of 
\begin{thm}[Key Lemma] \label{KeyLemma}
Let $X$ be a compact K\"ahler manifold with K\"ahler form $\omega$ and let
$\alpha \in H^{1,1}(X, \mathbb{R})$ be a pseudoeffective $(1,1)$ class. If
$\mathcal{F}$ and $\mathcal{G}$ are numerically trivial foliations then 
$\mathcal{F} \sqcup \mathcal{G}$ will be numerically trivial, too.
\end{thm}

\noindent The proof is divided in two main parts: First we show in 
Thm.~\ref{SingEst-thm} that we can neglect arbitrary (small neighborhoods of) 
analytic subsets when checking $(NT)_u$. This is done by applying Boucksom's
uniform bounds for the occuring integrals on generalized Lelong numbers.

\noindent The second part follows the inductive construction of 
$\mathcal{F} \sqcup \mathcal{G}$ on open subsets $U$ which do not intersect an 
arbitrarily small neighborhood $W$ of 
$\mathrm{Sing\ }\mathcal{F} \cup \mathrm{Sing\ }\mathcal{G} \cup 
 \mathrm{Sing\ }(\mathcal{F} \cap \mathcal{G})$ and the subset $Z$ analytic in 
$X-W$ on which the construction is not possible. If $z_1, \ldots, z_n$ are 
holomorphic coordinates on $U$ such that $\mathcal{F}_{|U}$ and 
$\mathcal{G}_{|U}$ are defined as in Prop.~\ref{FGcoord-prop} we prove first 
that the foliation $\mathcal{G}^\prime$ defined on $U$ by projection onto 
$z_{l+1}, \ldots , z_{n-m+l}$ is numerically trivial (Prop.~\ref{G'-prop}). 
Next 
we show that $\mathcal{F} + \mathcal{G}^\prime$ is numerically trivial in $U$
(Prop.~\ref{F+G-prop}; this is the first step of the inductive construction) 
and continue until we reach $\mathcal{F} \sqcup \mathcal{G}$.

\noindent We begin with
\begin{thm} \label{SingEst-thm}
Let $X$ be an $n$-dimensional compact K\"ahler manifold with K\"ahler form 
$\omega$ and let
$\alpha \in H^{1,1}(X, \mathbb{R})$ be a pseudoeffective $(1,1)$ class. Let
$U \subset X$ be an open subset with coordinates $z_1, \ldots, z_n$ and
$L = \{ z_1 = \ldots = z_l = 0 \} \subset U$ a linear subspace. Then for every
exhaustion $K_i \subset\subset K_{i+1} \subset\subset X$ of
$X - L = \bigcup K_i$ we have for all $1 \leq p \leq n-1$
\[ \lim_{\epsilon \downarrow 0} \sup_{T \in \alpha[-\epsilon\omega]} 
   \int_{(X - K_i) \cap U} (T_{ac} + \epsilon\omega)^p \wedge \omega^{n-p} 
   \stackrel{i \rightarrow \infty}{\longrightarrow} 0 \]
where the $T$'s run through all currents with analytic singularities 
representing $\alpha$ and $T \geq -\epsilon\omega$.
\end{thm}
\begin{proof}
Let us consider the generalized Lelong numbers of currents 
$(T_{ac} + \epsilon\omega)^p$ in points 
$x = (0, \ldots , 0, x_{l+1}, \ldots, x_n) \in L$ with respect to the 
plurisubharmonic weight
\[ \phi_x(z) = \log(\sum_{k=1}^{l} |z_i-x_i|^2+\sum_{k=l+1}^{n} |z_i-x_i|^q) \]
where $q$ is some integer $\geq n-l$. The advantage of this weight 
is 
that for a given $r > 0$ the number of subsets 
$\{z: |\phi_x(z)| \leq \log r\}$ necessary to cover $L$ is 
$\leq C \cdot r^{-\frac{2}{q} \dim L}$
for some constant $C > 0$ independent of $r$. Furthermore there are two 
constants
$C_1, C_2 > 0$
such that
\[ C_1\omega \leq \frac{i}{2}\partial\overline{\partial} e^{\phi_x} \leq 
   C_2\omega. \]

\noindent
By definition the generalized Lelong number 
$\nu((T_{ac} + \epsilon\omega)^p,\phi_x)$
is the decreasing limit for
$t \rightarrow -\infty$
of
\[ \begin{array}{rcl}
   \nu((T_{ac} + \epsilon\omega)^p,\phi_x,t) & = &  
   \int_{\phi_x(z)<t} (T_{ac} + \epsilon\omega)^p \wedge 
   (\partial\overline{\partial} \phi_x)^{n-p} \\
    & & \\
    & = & \frac{1}{(\pi e^{2t})^{n-p}} 
   \int_{\phi_x(z)<t} (T_{ac} + \epsilon\omega)^p \wedge 
   (\partial\overline{\partial} e^{\phi_x})^{n-p} 
   \end{array} \]
where the second equality follows from a formula proven in 
\cite[(2.13)]{Dem00}.

\noindent Set $t = \log r$. On the one hand, for
$r \leq r_0$
one has
\[ \nu((T_{ac} + \epsilon\omega)^p,\phi_x,\log r) \leq 
   \nu((T_{ac} + \epsilon\omega)^p,\phi_x,\log r_0) 
   \leq \frac{C_2}{(\pi r_0^2)^{n-p}} \int_X 
   (T_{ac} + \epsilon\omega)^p \wedge \omega^{n-p}. \]
But $T$ has analytic singularities. Hence using an idea of Boucksom 
\cite{Bou02}[3.1.12]  $\int_X (T_{ac} + \epsilon\omega)^p \wedge \omega^{n-p}$
is bounded from above by a constant depending only on the cohomology class 
$\alpha$ of $T$. On the other hand,
\[ (\pi r^2)^{n-p} \nu((T_{ac} + \epsilon\omega)^p,\phi_x ,\log r) \geq  C_1 
   \int_{e^{\phi_x(z)}<r} (T_{ac} + \epsilon\omega)^p \wedge \omega^{n-p}. \]
The claim follows from the choice of $l$ in the definition of $\phi_x$ and the 
upper bound on the number of level subsets covering $L$.
\end{proof}

\noindent From now on let $X$ be an $n$-dimensional compact K\"ahler manifold
with Kähler form $\omega$ and $\alpha \in H^{1,1}(X,\mathbb{R})$ a 
pseudo-effective $(1,1)$-class. We start the second part with the easiest 
configuration of two foliations $\mathcal{F}$ and $\mathcal{G}$ of $X$ in an 
open subset $U \subset X$:
\begin{prop} \label{F+G-prop}
Let $z_1, \ldots, z_n$ be coordinates on $U$ such that $\mathcal{F}$ is induced
by the projection on $z_{k+1}, \ldots, z_{n}$ and $\mathcal{G}$ by the 
projection on $z_{l+1}, \ldots , z_{n-m+l}$. If $\mathcal{F}$ and $\mathcal{G}$
are numerically trivial on $U$ the foliation $\mathcal{F} + \mathcal{G}$ 
induced by the projection on $z_{k+1}, \ldots , z_{n-m+l}$ is numerically 
trivial, too.
\end{prop}
\begin{proof}
By Theorem~\ref{NTCrit-thm} it is enough to show that constant $(n-p,n-p)$
test forms for $\mathcal{F} + \mathcal{G}$ are $\mathbb{C}$-linear combinations
of test forms for $\mathcal{F}$ and $\mathcal{G}$. But constant $(n-p,n-p)$
test forms for $\mathcal{F} + \mathcal{G}$ are $\mathbb{C}$-linear combinations
of decomposable  $(n-p,n-p)$ forms $dz_I \wedge d\overline{z}_J$ with
\[ \begin{array}{rcl}
   \left| I \cap \{k+1, \ldots, n-m+l\}\right| & >& n-m-k+l-p \mathrm{\ \ or}\\
   \left| J \cap \{k+1, \ldots, n-m+l\}\right| & > &  n-m-k+l-p. 
   \end{array} \]
Assume w.l.o.g. that the inequality for $I$ is satisfied and set
\[ q := \left| I \cap \{k+1, \ldots, n-m+l\}\right|. \]
If $dz_I \wedge d\overline{z}_J$ is not a test form for $\mathcal{F}$ and 
$\mathcal{G}$ we have
\[ \begin{array}{rcl}
   \left| I \cap \{k+1, \ldots, n\}\right| & \leq & n-k-p \mathrm{\ \ and}\\
   \left| I \cap \{l+1, \ldots, n-m+l\}\right| & \leq &  n-m-p. 
   \end{array} \]
This implies 
\[ \begin{array}{rcl}
   \left| I \cap \{l+1, \ldots, k\}\right| & \leq & n-m-p-q \mathrm{\ \ and}\\
   \left| I \cap \{n-m+l+1, \ldots, n\}\right| & > &  n-k-p-q. 
   \end{array} \]
Consequently
\[ \begin{array}{rcl}
   \left| I \right| & = & \left| I \cap \{1, \ldots, l\}\right| + 
                          \left| I \cap \{l+1, \ldots, k\}\right| + q + 
                          \left| I \cap \{n-m+l+1, \ldots, n\}\right| \\
        & \leq & l + n - m - p - q + q + n - k - p - q = 
                 n - p + (n - m - k + l - p) - q \\
        & < & n - p
   \end{array} \]
by the properties of $p$. This is a contradiction.
\end{proof}

\noindent Next we state a little fact that is useful later on again and again:
\begin{lem} \label{notF-lem}
Let $z_1, \ldots, z_n$ be coordinates on $U$ such that $\mathcal{F}$ is induced
by the projection on $z_{k+1}, \ldots, z_{n}$ and $\mathcal{G}$ by the 
projection on $z_{l+1}, \ldots , z_{n-m+l}$. If $dz_I \wedge d\overline{z}_J$
is not a (constant) $(n-p,n-p)$ test form for $\mathcal{F}$ we have
\[ \{l+1, \ldots, k\} \subset I, J. \]
\end{lem}
\begin{proof}
If $dz_I \wedge d\overline{z}_J$ is not a test form for $\mathcal{F}$ it 
follows that
\[ \left| I \cap \{k+1, \ldots, n\}\right|, 
   \left| J \cap \{k+1, \ldots, n\}\right| \leq n-k-p. \]
But then $|I| = |J| = n-p$ implies
\[ \left| I \cap \{1, \ldots, k\}\right| = k. \]
\end{proof}

\noindent From now on fix coordinates $z_1, \ldots, z_n$ on $U$ such that the 
foliations $\mathcal{F}$ and $\mathcal{G}$ are described as in 
Proposition~\ref{FGcoord-prop}. 

\noindent The main idea for proving the numerical triviality of 
$\mathcal{G}^\prime$ is to compare the evaluation of test forms on fibers of 
the projection $\pi : U \rightarrow \mathbb{C}^{k-l}$ onto 
$z_{l+1}, \ldots , z_{k}$. Intuitively the numerical triviality should imply 
that the difference of these values vanishes. Since we always compute integrals
on $U$ we first need a comparison lemma for \textit{nearby} fibers:
\begin{lem} \label{ApproxPush-lem}
Suppose that $0 \in U \subset \mathbb{C}^n$. Let 
\[ \Phi : U \rightarrow \mathbb{C}^n, (z_1, \ldots , z_n) \mapsto
          (z_1, \ldots , z_k, z_{k+1}^\prime = \Phi_{k+1}(z), \ldots, 
                              z_n^\prime = \Phi_n(z)) \]
be a coordinate transformation such that
\[ \Phi_{k+j}(z)_{|U_0} = z_{k+j}\ \mathrm{on\ } 
   U_0 = \{z \in U : z_{l+1} = \ldots = z_{k} = 0 \}. \]
Set $U_\delta := \pi^{-1}(B_\delta(0))$ and note that 
$\Phi(U_\delta)=U_\delta$. Then for every real $(n-p,n-p)$-form $u$ with 
$p \leq n-k+l$ which is not a test form for $\mathcal{F}$ there exists a 
constant $C > 0$ independent of $\delta$ such
that for every $\delta$ small enough the inequality of $(n-p,n-p)$-forms
\[ -C\delta\omega^{n-p-k+l}_{|\pi^{-1}(x)} < (\Phi_\ast(u) - u)_{|\pi^{-1}(x)} 
   < C\delta\omega^{n-p-k+l}_{|\pi^{-1}(x)} \]
is true on all fibers $\pi^{-1}(x)$, $x \in B_\delta(0)$. 
\end{lem}
\begin{proof}
We can replace $\omega$ by the standard $(1,1)$-form 
$\sum dz_i \wedge d\overline{z}_i$ on $\mathbb{C}^n$. Let $f_{IJ}$ be the 
coefficient functions of $(\Phi_\ast(u) - u)_{|\pi^{-1}(x)}$ w.r.t. the base 
$dz_I \wedge d\overline{z}_J$. Then $(\Phi_\ast(u) - u)_{|\pi^{-1}(0)} = 0$
implies $f_{IJ|\pi^{-1}(0)} \equiv 0$.
Let $v \in \mathbb{C}^n$ be a direction vector with 0 entries in all 
coordinates but $z_{l+1}, \ldots, z_k$ and $\parallel\!\! v \!\!\parallel = 1$.
The mean value theorem gives us for all $x_0 \in \pi^{-1}(0)$
\[ \parallel\!\! f_{IJ}(x_0 + t_0v) \!\!\parallel = 
   | \frac{d}{dt} f_{IJ}(x_0 +tv) \cdot t_0v | = 
   |t_0| \cdot \parallel\!\! D_{(l+1, \ldots, k)} f_{IJ}(x_0 +tv) \cdot v 
           \!\!\parallel \]
where the matrices $D_{(l+1, \ldots, k)} f_{IJ}$ collect the partial 
derivatives
w.r.t. $x_{l+1}, y_{l+1}, \ldots, x_k, y_k$ ($z_j = x_j + iy_j$). These 
matrices can be considered as continuous families of linear maps and hence 
their norms are
bounded from above by a constant $C^\prime$. Consequently
\[ \parallel\!\! (\Phi_\ast(u) - u)_{|\pi^{-1}(x)} \!\!\parallel_{\sup} 
   \ \leq\  C^\prime \delta \]
and the claim follows.
\end{proof}

\begin{rem}
$\omega_{|\pi^{-1}(x)}$ is the usual restricted form on the submanifold 
$\pi^{-1}(x)$. But $(\Phi_\ast(u) - u)_{|\pi^{-1}(x)}$ must be defined as the
$(n-p-k+l,n-p-k+l)$-form obtained from $\Phi_\ast(u) - u$ by replacing 
$dz_I \wedge d\overline{z}_J$ with 
$dz_{I-\{l+1, \ldots, k\}} \wedge d\overline{z}_{J-\{l+1, \ldots, k\}}$. This 
makes sense by 
Lemma~\ref{notF-lem} because $u$ is not a test form for $\mathcal{F}$.
\end{rem}

\begin{prop} \label{nearbyfiber-prop}
Let $x$ be any point in $U$.
If $u$ is a real constant $(n-p,n-p)$ test form for $\mathcal{G}^\prime$ in $U$
but not for $\mathcal{F}$ and $p \leq n-k+l$ then
\[ \lim_{\delta \rightarrow 0} \frac{1}{\mathrm{Vol}(U_\delta)} 
   \lim_{\epsilon \downarrow 0} \sup_{T \in \alpha[-\epsilon\omega]}
   \int_{U_\delta} \left| (T_{ac}+\epsilon\omega)^p \wedge u \right| = 0 \]
where $U_\delta \subset U$  denotes the open subset 
$\pi^{-1}(B_\delta(\pi(x)))$.
\end{prop}
\begin{proof}
Using the notation of Proposition~\ref{FGcoord-prop} the map 
\[ \Phi : U \rightarrow \mathbb{C}^n , (z_1, \ldots , z_n) \mapsto 
   (z_1, \ldots , z_k, g_{m+k-l+1}(z), \ldots , g_n(z), 
    z_{n-m+l+1}, \ldots , z_n) \]
describes a coordinate transformation such that $u^\prime := \Phi_\ast(u)$
is a real constant test form for $\mathcal{G}$. Hence
\[ (\ast) \rule{2cm}{0cm} 
    \lim_{\epsilon \downarrow 0} \sup_{T \in \alpha[-\epsilon\omega]}
       \int_U \left| (T_{ac}+\epsilon\omega)^p \wedge u^\prime \right| = 0. 
           \rule{5cm}{0cm} \]
We replace $U$ by $U_\delta$ and want to compare the growth of this limit
with that of $\mathrm{Vol}(U_\delta)$. To this purpose it is enough to look at
sequences of currents $T_k \in \alpha[-\epsilon_k\omega]$ with analytic 
singularities and 
$\epsilon_k \stackrel{k \rightarrow \infty}{\longrightarrow} 0$. 

\noindent Since $u$ is not a test form for $\mathcal{F}$ Lemma~\ref{notF-lem}
tells us that $u$ only contains forms $dz_I \wedge d\overline{z}_J$ with
$\{l+1, \ldots, k\} \subset I,J$. The same is true for 
$u^\prime = \Phi_\ast(u)$ in the new coordinates. Hence using Fubini's theorem
\[ \int_U \left| ((T_k)_{ac}+\epsilon_k\omega)^p \wedge u^\prime \right| = 
   \int_{\pi(U)} d\lambda_{\pi(U)}(x) 
   \left( \int_{\pi^{-1}(x)} 
   \left| ((T_k)_{ac}+\epsilon_k\omega)_{|\pi^{-1}(x)}^p \wedge 
                                         u^{\prime\prime} \right| \right) \]
for the $(n-k+l-p, n-k+l-p)$-form $u^{\prime\prime}$ obtained from $u$ as 
described in the remark above.

\noindent Define $L^1$-functions
\[ f_k : \pi(U) \rightarrow \mathbb{R}^+, \ 
   x \mapsto \int_{\pi^{-1}(x)} 
   \left| ((T_k)_{ac}+\epsilon_k\omega)_{|\pi^{-1}(x)}^p \wedge 
                                         u^{\prime\prime} \right|. \]
The $f_k$'s tend to $0$ in $L^1$-norm, by $(\ast)$. Convoluting the $f_k$'s
with $\rho_\delta := \frac{1}{\mathrm{Vol}(U_\delta)} \chi_{U_\delta}$ we get
\[ \frac{1}{\mathrm{Vol}(U_\delta)} \cdot \lim_{k \rightarrow \infty}
   \int_{U_\delta} \left| ((T_k)_{ac}+\epsilon_k\omega)^p 
                                          \wedge u^\prime \right| = 
   \lim_{k \rightarrow \infty} (\rho_\delta \ast f_k)(0) = 0 \] 
because convolution with characteristic functions of open subsets improves
$L^1$-convergence in sup-norm convergence. We conclude
\begin{eqnarray*}
\lefteqn{\frac{1}{\mathrm{Vol}(U_\delta)} 
   \lim_{\epsilon \downarrow 0} \sup_T
   \int_{U_\delta} \left| (T_{ac}+\epsilon\omega)^p \wedge u \right|\ \leq}\\
     & \leq & \frac{1}{\mathrm{Vol}(U_\delta)} \left(
       \lim_{\epsilon \downarrow 0} \sup_T
       \int_{U_\delta} \left| (T_{ac}+\epsilon\omega)^p \wedge 
                                                   (u-u^\prime) \right| +
       \lim_{\epsilon \downarrow 0} \sup_T
       \int_{U_\delta} \left| (T_{ac}+\epsilon\omega)^p \wedge u^\prime \right|
       \right)\\
 & \leq & \frac{1}{\mathrm{Vol}(U_\delta)} \lim_{\epsilon \downarrow 0} \sup_T
   \int_{B_\delta(0)} d\lambda_{B_\delta(0)}(x) \left( \int_{\pi^{-1}(x)} 
   \left| (T_{ac}+\epsilon\omega)_{|\pi^{-1}(x)}^p \wedge 
                                      (u-u^\prime)_{|\pi^{-1}(x)} \right| 
        \right)
\end{eqnarray*}

\noindent Using Lemma~\ref{ApproxPush-lem} we continue this inequality chain
with
\begin{eqnarray*}
 & \leq & \frac{C}{\mathrm{Vol}(U_\delta)} \lim_{\epsilon \downarrow 0} \sup_T
   \int_{B_\delta(0)} d\lambda_{B_\delta(0)}(x) \left( \int_{\pi^{-1}(x)}
   (T_{ac}+\epsilon\omega)_{|\pi^{-1}(x)}^p \wedge 
                                      \delta\omega^{n-p-k+l}_{|\pi^{-1}(x)} 
   \right)\\
 & \leq & \frac{C \cdot \delta}{\mathrm{Vol}(U_\delta)}
          \lim_{\epsilon \downarrow 0} \sup_T \int_{U_\delta}
          (T_{ac}+\epsilon\omega)^p \wedge \omega^{n-p} \\
 & \leq & \frac{C^\prime \cdot \delta}{\delta^{2(k-l)}}
          \lim_{\epsilon \downarrow 0} \sup_T \delta^{2(n-p)} \cdot 
   \nu((T_{ac} + \epsilon\omega)^p,\log |z|,\log \delta) \\
 & \leq & C^\prime \cdot \delta \cdot \lim_{\epsilon \downarrow 0} \sup_T
          \int_X (T_{ac}+\epsilon\omega)^p \wedge \omega^{n-p}
\end{eqnarray*}
by the estimates in the proof of Theorem~\ref{SingEst-thm} and $p \leq n-k+l$.
This last term tends to $0$ if $\delta \rightarrow 0$ because of the uniform 
bounds of Boucksom.
\end{proof}

\noindent Now we compare different fibers of $\pi$:
\begin{prop} \label{fibercomp-prop}
Assume that $x = 0$ and let $u$ be again a real $(n-p,n-p)$ form, 
$p \leq n-k+l$, which is not a test form for $\mathcal{F}$. Let $\overline{y}$ 
be another point in $\pi(U) \subset \mathbb{C}^{k-l}$ and 
$U_{\delta,\overline{y}} := \pi^{-1}(B_\delta(\overline{y}))$. Then
\[ \int_{U_{\delta,\overline{y}}} 
   \left| (T_{ac}+\epsilon\omega)^p \wedge u \right| =
   \int_{U_\delta} \left| (T_{ac}+\epsilon\omega)^p \wedge u \right| . \]
\end{prop}
\begin{proof}
Let
$\overline{x} = (z_{l+1}, \ldots, z_k) = 
   (x_{l+1}, y_{l+1}, \ldots , x_k, y_k) \in B_\delta(0)$. Then 
$\pi^{-1}(\overline{x})$ and $\pi^{-1}(\overline{x}+\overline{y})$ may be seen 
as part of the boundary of a ``cylinder'' $Z(\overline{x})$ obtained by 
connecting all pairs of points $x^\prime \in \pi^{-1}(\overline{x})$ and
$x^\prime + (0,\ldots,0,\overline{y},0,\ldots,0) \in 
 \pi^{-1}(\overline{x}+\overline{y})$ with a real line.

\noindent Since $T_{ac}$ only has analytic singularities of codimension 
$\geq 2$ the restriction of $T_{ac}$ to these real lines is a smooth form for
almost all $x^\prime \in \pi^{-1}(\overline{x})$.

\noindent To prove the proposition it is enough to look at $u$'s which are 
decomposable forms $dx_{I_0} \wedge dy_{J_0}$ in the real coordinates 
$x_1,\ldots,x_n,y_1,\ldots,y_n$. Then
\begin{eqnarray*} 
   (T_{ac}+\epsilon\omega)^p \wedge dx_{I_0} \wedge dy_{J_0} & = &  
   \left(\sum_{|I|=|J|=p} T_{IJ} dx_I \wedge dy_J\right) 
                                             \wedge dx_{I_0} \wedge dy_{J_0} \\
   & = & T_{I_0^\prime J_0^\prime} dx_{I_0^\prime} \wedge dy_{J_0^\prime} 
                                             \wedge dx_{I_0} \wedge dy_{J_0} 
\end{eqnarray*}
where $I_0^\prime \cup I_0 = J_0^\prime \cup J_0 = \{1,\ldots,n\}$ and 
$T_{I_0^\prime J_0^\prime}$ is a real $L^1_{loc}$-function.

\noindent Let $\overline{u}$ be the decomposable form for which 
\[ u = \overline{u} \wedge x_{l+1} \wedge y_{l+1} \wedge \ldots 
                     \wedge x_k \wedge y_k. \]
The form $\overline{u}$ exists because $u$ is not a test form for 
$\mathcal{F}$. Then
\[ (T_{ac}+\epsilon\omega)^p \wedge \overline{u} = T_{I_0^\prime J_0^\prime}
   dx_I \wedge dy_J,\ I = J = \{1,\ldots,l,k+1,\ldots,n\}. \] 
Since $T_{ac}$, $\omega$ and $\overline{u}$ (as a constant form) are closed 
forms resp. currents we get
\[ d \left[ (T_{ac}+\epsilon\omega)^p \wedge \overline{u} \right] = 0. \]
If we restrict $(T_{ac}+\epsilon\omega)^p \wedge \overline{u}$ to the cylinder
$Z(\overline{x})$ the current remains closed. But in $Z(\overline{x})$ this 
just means that the derivation in direction of 
$(0,\ldots,0,\overline{y},0,\ldots,0)$ vanishes. This implies that 
$T_{I_0^\prime J_0^\prime}$ remains constant on the real lines in 
$Z(\overline{x})$ where $T_{ac}$ is smooth. Consequently integrating
$\left| (T_{ac}+\epsilon\omega)^p \wedge \overline{u} \right|$ on the top
$\pi^{-1}(\overline{x}+\overline{y})$ and the bottom $\pi^{-1}(\overline{x})$
of the cylinder $Z(\overline{x})$ gives the same result. We finally calculate:
\begin{eqnarray*}
\lefteqn{\int_{U_{\delta,\overline{y}}} 
   \left| (T_{ac}+\epsilon\omega)^p \wedge u \right| -
   \int_{U_\delta} \left| (T_{ac}+\epsilon\omega)^p \wedge u \right| = } \\
   & = & \int_{B_\delta(0)} d\overline{x} \left( 
         \int_{\pi^{-1}(\overline{x}+\overline{y})} 
         \left| (T_{ac}+\epsilon\omega)^p \wedge \overline{u} \right| -
         \int_{\pi^{-1}(\overline{x})} 
         \left| (T_{ac}+\epsilon\omega)^p \wedge \overline{u} \right|
         \right) = 0. 
\end{eqnarray*}
\end{proof}

\noindent Now we can show
\begin{prop} \label{G'-prop}
If $\mathcal{F}$ and $\mathcal{G}$ are numerically trivial in $U$ w.r.t. 
$\alpha$ then the foliation $\mathcal{G}^\prime$ will be numerically trivial in
$U$ w.r.t. $\alpha$, too.
\end{prop}
\begin{proof}
We have to show that
\[ \lim_{\epsilon \downarrow 0} \sup_{T \in \alpha[-\epsilon\omega]}
   \int_U \left| (T_{ac}+\epsilon\omega)^p \wedge u \right| = 0 \]
for every constant $(n-p,n-p)$-test form for  $\mathcal{G}^\prime$. Since
$\mathcal{F}$ is numerically trivial in $U$ we only need to check forms which 
are not test forms for $\mathcal{F}$. To this purpose we cover $U$ with open
subsets $U_\delta^{(i)} := U_{\delta,\overline{y}_i}$ as in the previous 
proposition. There exists a constant $C > 0$ independent of $\delta$ such that
we only need $\frac{C}{\mathrm{Vol(U_\delta)}}$ of these covering sets for 
every $\delta > 0$. Now we calculate:
\begin{eqnarray*}
\lim_{\epsilon \downarrow 0} \sup_{T \in \alpha[-\epsilon\omega]}
   \int_U \left| (T_{ac}+\epsilon\omega)^p \wedge u \right| & \leq &
\lim_{\epsilon \downarrow 0} \sup_{T \in \alpha[-\epsilon\omega]} \sum_i
   \int_{U_\delta^{(i)}} \left| (T_{ac}+\epsilon\omega)^p \wedge u \right| \\
 & \leq & \sum_i 
   \lim_{\epsilon \downarrow 0} \sup_{T \in \alpha[-\epsilon\omega]}
   \int_{U_\delta^{(i)}} \left| (T_{ac}+\epsilon\omega)^p \wedge u \right|
\end{eqnarray*}
Applying Proposition~\ref{fibercomp-prop} the last term equals
\[ \sum_i \lim_{\epsilon \downarrow 0} \sup_T
   \int_{U_\delta} \left| (T_{ac}+\epsilon\omega)^p \wedge u \right| \leq
   \frac{C}{\mathrm{Vol(U_\delta)}} 
   \sum_i \lim_{\epsilon \downarrow 0} \sup_T
   \int_{U_\delta} \left| (T_{ac}+\epsilon\omega)^p \wedge u \right|
\]
which tends to $0$ for $\delta \rightarrow 0$ by Prop.~\ref{nearbyfiber-prop}.
\end{proof}

\noindent Applying Prop.~\ref{G'-prop} and Prop.~\ref{F+G-prop} in each step of
the inductive construction of $\mathcal{F} \cup \mathcal{G}$ we finish the 
proof of the Key Lemma~\ref{KeyLemma}. \hfill $\Box$

\noindent The Key Lemma allows to construct a maximal foliation numerically 
trivial w.r.t. $\alpha$ which will be called \textit{the} numerically trivial
foliation w.r.t. $\alpha$. We are now able to prove Thm.~\ref{AppBd-thm}
from the introduction:
\begin{thm}
Let $X$ be an $n$-dimensional compact K\"ahler manifold with K\"ahler form 
$\omega$ and let $\alpha \in H^{1,1}(X, \mathbb{R})$ be a pseudoeffective 
$(1,1)$-class with numerical dimension $\nu(\alpha)$. Let $\mathcal{F}$ be the 
numerically trivial foliation w.r.t. $\alpha$. Then
\[ \mathrm{rk}(\mathcal{F}) \leq n - \nu(\alpha). \]
\end{thm}
\begin{proof}
Set $k := \mathrm{rk}(\mathcal{F})$. If $n-k < \nu(\alpha)$ every 
$(n-\nu(\alpha),n-\nu(\alpha))$-form with compact support in 
$X - \mathrm{Sing}(\mathcal{F})$ is a test form for $\mathcal{F}$. In 
particular, for an arbitrarily small compact subset 
$K \supset \mathrm{Sing}(\mathcal{F})$ we have
\[ \lim_{\epsilon \downarrow 0} \sup_{T \in \alpha[-\epsilon\omega]}
   \int_{X-K} (T_{ac}+\epsilon\omega)^{\nu(\alpha)} \wedge 
                                      \omega^{n-\nu(\alpha)} = 0. \]
By Thm.~\ref{SingEst-thm}
\[ \lim_{\epsilon \downarrow 0} \sup_{T \in \alpha[-\epsilon\omega]}
   \int_{K} (T_{ac}+\epsilon\omega)^{\nu(\alpha)} \wedge 
                                      \omega^{n-\nu(\alpha)} \rightarrow 0 \]
uniformly with the volume of $K$. Consequently
\[ \lim_{\epsilon \downarrow 0} \sup_{T \in \alpha[-\epsilon\omega]}
   \int_X (T_{ac}+\epsilon\omega)^{\nu(\alpha)} \wedge 
                                      \omega^{n-\nu(\alpha)} = 0. \]
That contradicts the definition of $\nu(\alpha)$.
\end{proof}

\section{The Transversality Lemma}

\noindent It is difficult to determine the numerically trivial foliation of a 
pseudo-effective $(1,1)$-class. Sometimes the following theorem helps:
\begin{thm}[Transversality Lemma] \label{transvers-lem}
Let $X$ be an $n$-dimensional compact K\"ahler manifold with K\"ahler form 
$\omega$, and let $\alpha \in H^{1,1}(X,\mathbb{R})$ be a pseudo-effective 
$(1,1)$-class with numerical dimension $\nu(\alpha)$.

\noindent Let $\mathcal{F}$ be a foliation of rank $k=n-\nu(\alpha)$ and 
$\{U\}$ a covering of $X-\mathrm{Sing}(\mathcal{F})$ by open subsets $U$ with
coordinates $z_1, \ldots , z_n$ such that the leaves of $\mathcal{F}$ are the
fibers of the projection on $z_{k+1}, \ldots, z_n$.

\noindent Now suppose that for all $\epsilon>0$ and all open subsets $U$ of the
covering there exists a constant $\delta_U>0$ and a current 
$T_{\epsilon,U} \in \alpha[-\epsilon\omega]$ such that 
\[ (T_{\epsilon,U} + \epsilon\omega)_{|U} \geq 
   \delta_U \cdot \omega_{\mathcal{F},U} := 
   \delta_U \cdot \sum_{j=k+1}^n dz_j \wedge d\overline{z}_j. \]
Then $\mathcal{F}$ is the numerically trivial foliation of $\alpha$.
\end{thm}
\begin{proof}
By Thm.~\ref{NTCrit-thm} it is enough to show on every $U$ of the covering that
\[ (\ast)  \rule{2cm}{0cm}  
   \lim_{\epsilon \downarrow 0} \sup_{T \in \alpha[-\epsilon\omega]}
   \int_U \left| (T_{ac}+\epsilon\omega)^p \wedge 
                                  dz_I \wedge d\overline{z}_J \right| = 0 
   \rule{5cm}{0cm} \]
for every $1 \leq p \leq n-1$ and every constant $(n-p,n-p)$ test form 
$dz_I \wedge d\overline{z}_J$ for $\mathcal{F}$ on $U$. We start with proving
$(\ast)$ for test forms $dz_I \wedge d\overline{z}_I$ with
\[ m := |I \cap \{k+1, \ldots, n\}| > n-k-p. \]
From $n-k = \nu(\alpha)$ and the definition of the numerical dimension we 
conclude that
\[ \lim_{\epsilon \downarrow 0} \sup_{T \in \alpha[-\epsilon\omega]}
   \int_X (T_{ac}+\epsilon\omega)^p \wedge 
          (T_{\epsilon,U,\mathrm{ac}}+\epsilon\omega)^m \wedge 
          \omega^{n-p-m} = 0.  \]
The inequality $(T_{\epsilon,U} + \epsilon\omega)_{|U} \geq 
   \delta_U \cdot \omega_{\mathcal{F},U}$ implies 
\[ \lim_{\epsilon \downarrow 0} \sup_{T \in \alpha[-\epsilon\omega]}
   \int_X (T_{ac}+\epsilon\omega)^p \wedge 
          \omega_{\mathcal{F},U}^m \wedge \omega^{n-p-m} = 0.\]
But there exists a constant $C>0$ such that
\[ \omega_{\mathcal{F},U}^m \wedge \omega^{n-p-m} > 
   \pm C \cdot i^{n-p} dz_I \wedge d\overline{z}_I \]
hence $(\ast)$ for $dz_I \wedge d\overline{z}_I$.

\noindent Next we look at general test forms $dz_I \wedge d\overline{z}_J$. We
can assume w.l.o.g. that
\[ |I \cap \{k+1,\ldots,n\}| > n-k-p,\ |J \cap \{k+1,\ldots,n\}| \geq n-k-p.\] 
Since
\[ (T_{ac}+\epsilon\omega)^p = 
   \sum_{|I^\prime|=|J^\prime|=p} T_{I^\prime J^\prime} 
   dz_{I^\prime} \wedge d\overline{z}_{J^\prime} \]
is a smooth semipositive form outside an analytic subset we have
\[ |T_{I^\prime J^\prime}| \leq |T_{I^\prime I^\prime}|^{\frac{1}{2}} \cdot 
                                |T_{J^\prime J^\prime}|^{\frac{1}{2}} \]
almost everywhere by the Cauchy-Schwarz inequality. Hence for 
$I^\prime = \{1,\ldots,n\} - I$, $J^\prime = \{1,\ldots,n\} - J$ we get
\begin{eqnarray*}
\lefteqn{\int_U \left| (T_{ac}+\epsilon\omega)^p \wedge 
                                  dz_I \wedge d\overline{z}_J \right|  = 
\int_U |T_{I^\prime J^\prime}| idz_1 \wedge d\overline{z}_1 \wedge \ldots 
                        \wedge idz_n \wedge d\overline{z}_n} \\
 & \leq & \left( \int_U |T_{I^\prime I^\prime}|dV\right)^{\frac{1}{2}} \cdot
          \left( \int_U |T_{J^\prime J^\prime}|dV\right)^{\frac{1}{2}} \\
 & = & \left( \int_U \left| (T_{ac}+\epsilon\omega)^p \wedge 
                dz_I \wedge d\overline{z}_I \right| \right)^{\frac{1}{2}} \cdot
       \left( \int_U \left| (T_{ac}+\epsilon\omega)^p \wedge 
                dz_J \wedge d\overline{z}_J \right| \right)^{\frac{1}{2}} \\
 & \leq & \left( \int_U \left| (T_{ac}+\epsilon\omega)^p \wedge 
                dz_I \wedge d\overline{z}_I \right| \right)^{\frac{1}{2}} \cdot
       \left( \int_U (T_{ac}+\epsilon\omega)^p \wedge 
                \omega^{n-p} \right)^{\frac{1}{2}}
\end{eqnarray*}
where the first inequality is a consequence of the H\"older inequality. By 
Boucksom's uniform estimates the second integral is uniformly bounded from 
above and the first factor tends to $0$ when $\epsilon \rightarrow 0$ by what 
we have shown before. 

\noindent Finally we must exclude the possibility that some foliation 
$\mathcal{F}^\prime \supset \mathcal{F}$ different from $\mathcal{F}$ is 
numerically trivial w.r.t. $\alpha$. So let 
$\mathrm{rk}(\mathcal{F}^\prime) =: k^\prime > k$ and choose an open subset 
$U \subset X - (\mathrm{Sing}(\mathcal{F}) \cup 
                \mathrm{Sing}(\mathcal{F}^\prime))$ with coordinates 
$z_1, \ldots, z_n$ such that the leaves of $\mathcal{F}^\prime$ are the fibers
of the projection on $z_{k^\prime+1}, \ldots, z_n$ and the leaves of 
$\mathcal{F}$ the fibers of the projection on $z_{k+1}, \ldots, z_n$. This is
possible because of Prop.~\ref{GsubF-prop}.

\noindent Now consider the $(n-1,n-1)$ form 
$\eta = \pm i^{n-1} dz_I \wedge d\overline{z}_I$ given by
\[ I = \{1, \ldots, k^\prime-1, k^\prime+1, \ldots, n\}. \]
$\eta$ is a test form for $\mathcal{F}^\prime$ since
$|I \cap \{k^\prime+1, \ldots, n\} = n-k^\prime > n-k^\prime-1$. But there 
exists a constant $C>0$ such that
\[ \int_U \left| (T_{ac}+\epsilon\omega) \wedge \eta \right| \geq \delta_u 
   \int_U \left| \omega_{\mathcal{F},U} \wedge \eta \right| \geq C \cdot 
   \int_U \omega^n > 0. \]
Hence $\mathcal{F}^\prime$ is not numerically trivial w.r.t. $\alpha$.
\end{proof} 

\noindent As a first application we use that the foliations 
constructed on the surface examples 4.1 and 4.2 in \cite{Eck04b} satisfy the 
conditions of the Transversality Lemma and conclude that they are still the
numerically trivial foliations.

\noindent More general the definition of numerical triviality in \cite{Eck04b}
implies Def.~\ref{NTFol-Def} on surfaces.

\noindent The Transversality Lemma gives also a very simple proof for
\begin{thm} \label{KodIitfib=NTFol-thm}
Let
$X$
be a K\"ahler manifold and 
$L$
a pseudo-effective line bundle on
$X$.
Suppose that the Kodaira-Iitaka dimension
$\kappa(X,L)$
equals the numerical dimension 
$\nu(X,L)$
of 
$L$.
Then the numerical trivial foliation of
$L$
is the Kodaira-Iitaka fibration of
$L$.
\end{thm}
\begin{proof}
Choose $m \gg 0$ such that the linear system $|mL|$ defines the Kodaira-Iitaka 
fibration $f: X \dasharrow Y \subset \mathbb{P}^N$. Then the positive curvature
current of the singular metric $h_{|mL|}$ satisfies the conditions in the 
Transversality Lemma because it is the pull back of the Fubini-Study metric on
$\mathcal{O}(1)$ over $Y$.
\end{proof}

\section{Variants and derived constructions}

\noindent Obviously numerical triviality may also be defined with $T$ running 
through subsets $\mathcal{C}(\epsilon) \subset \alpha[-\epsilon\omega]$ of
$(1,1)$-currents with analytic singularities such that
\[ \mathcal{C}(\epsilon^\prime) \subset \mathcal{C}(\epsilon) \]
for $0 < \epsilon^\prime \leq \epsilon$. The uniform boundedness of the mass of
$T^p_{ac}$ for arbitrary currents $T \in \alpha[-\epsilon\omega]$ 
(\cite[Thm. 3.1.10]{Bou02}) shows that we can even use a fixed positive current
$T \in  \alpha[0]$ with \textit{arbitrary} singularities. In this case $(NT)_u$
reduces to 
\[ \int_X |T^p_{ac} \wedge u| = 0 \]
for all test forms $u$ of a foliation $\mathcal{F}$ and 
Lemma~\ref{pullback-lem} implies that locally $T_{ac}$ is the pull back of a 
current $S$ on the base of the projection locally defining $\mathcal{F}$. In 
particular this numerically trivial foliation w.r.t. $T$ equals the one defined
in \cite[Def.2.9]{Eck04b}.

\noindent Furthermore let $(T_k)_{k \in \mathbb{N}}$ be a sequence of currents
$T_k \in \alpha[-\epsilon_k\omega]$ with analytic singularities such that 
$\epsilon_k \stackrel{k \rightarrow \infty}{\longrightarrow} 0$ and 
$T_{k,ac} \rightarrow T_{ac}$ almost everywhere. The Fatou lemma shows that
\[ \int_X |T^p_{ac} \wedge u| \leq \liminf_{k \rightarrow \infty}
   \int_X |(T_{k,ac}+\epsilon_k\omega)^p \wedge u| \]
hence the numerically trivial foliation w.r.t. $(T_k)_{k \in \mathbb{N}}$ is 
contained in the numerically trivial foliation w.r.t. $T$. The same is true for
the numerically trivial foliation of $\alpha$.

\noindent The inclusion can be strict as the surface example \cite[4.1]{Eck04b}
shows: The smooth currents $T_k \in c_1(L)[-\epsilon_k\omega]$ constructed 
there define a numerically trivial foliation with $1$-dimensional leaves. By 
weak compactness of almost positive $(1,1)$-currents there exists a subsequence
of the $T_k$ weakly converging to a positive $(1,1)$-current in $c_1(L)$. But
the only positive current in $c_1(L)$ is the integration current of a divisor
whose absolutely continuous part is $0$.

\noindent Considering the Iitaka fibration of a line bundle $L$ with Kodaira 
dimension $\kappa(L) \geq 0$ it was proven in \cite[2.4]{Eck04b} that this 
fibration is the numerically trivial foliation w.r.t. the curvature current
$\Theta_{|mL|}$ of the positive metric $h_{mL}$  
defined by all sections of $|mL|$ (for an appropriate $m \gg 0$). Hence the 
numerically trivial foliation of
$c_1(L)$ is contained in the Iitaka fibration of $L$.

\noindent In general it is not true that numerically trivial foliations are 
(rational) fibrations, see the surface examples in \cite{Eck04b}. This 
motivates the following
\begin{Def}
Let $X$ be a compact K\"ahler manifold and $\alpha \in H^{1,1}(X,\mathbb{R})$ a
pseudo-effective $(1,1)$-class. Let $\mathcal{F}$ be the numerically trivial 
foliation of $\alpha$. Then the maximal meromorphic map $f: X \dasharrow Y$ 
such that the induced foliation is contained in $\mathcal{F}$ is called the 
\textbf{pseudo-effective fibration} of $\alpha$.
\end{Def}

\noindent Note that
\begin{itemize}
\item[(1)] 
Prop.~\ref{fibjoin-prop} shows that the definition makes sense: There is a 
maximal fibration contained in a foliation.
\item[(2)] 
The same definition for numerically trivial foliations w.r.t. a single 
positive current leads to Tsuji's numerically trivial fibrations (see 
\cite{Eck04a}) by \cite[Prop.2.12]{Eck04b}.
\end{itemize}

\noindent In the projective case we have an algebraic characterization of the
pseudo-effective fibration. It uses Boucksom's construction of a divisorial 
Zariski decomposition (\cite[Ch.2]{Bou02}): For every pseudo-effective line 
bundle $L$ there is a real $(1,1)$-class $Z(L)$ nef in codimension $1$ and a 
negative part $N(L) = \sum_E \nu(L,E) E$ where $\nu(L,E) \neq 0$ for only 
finitely many exceptional divisors $E \in X$ such that
\begin{itemize}
\item[(i)] $c_1(L) = Z(L) + N(L)$ and
\item[(ii)] the natural map 
$H^0(X, \lfloor kZ(L) \rfloor) \rightarrow H^0(X,kL)$ is an isomorphism for all
$k \in \mathbb{N}$.
\end{itemize}

\begin{prop} \label{pseff-fib-prop}
Let $X$ be a complex projective manifold and $L$ a pseudo-effective line 
bundle. Then there exists a covering family $(C_t)_{t \in T}$ of curves in $X$
such that 
\begin{itemize}
\item[(i)] $C_t.(L-N(L)) = 0$ and
\item[(ii)] the pseudo-effective fibration $f$ of $L$ resp. $\alpha := c_1(L)$
is the generic reduction map w.r.t. $(C_t)_{t \in T}$
\end{itemize}

\noindent Furthermore for every irreducible curve through a general point 
$x \in X$ not lying in a fiber of $f$ we have
\[ L.C > 0. \]
\end{prop}
\begin{proof}
Since foliations are uniquely determined by their restriction to Zariski-open 
subsets (just saturate) and numerical triviality can be checked outside 
analytic subsets (by Theorem~\ref{SingEst-thm}) numerically trivial foliations 
behave well under certain types of holomorphic maps: Let $X$ be  a compact 
K\"ahler manifold and $\alpha \in H^{1,1}(X,\mathbb{R})$ a pseudo-effective 
$(1,1)$-class.
\begin{itemize}
\item[(a)]
Let $\pi : \widehat{X} \rightarrow X$ be a modification of compact K\"ahler 
manifolds. Then the numerically trivial foliation of $\pi^\ast \alpha$ on 
$\widehat{X}$ is the pull back of the numerically trivial foliation of 
$\alpha$ on $X$, and vice versa.
\item[(b)]
The same is true for branched coverings $\pi: Y \rightarrow X$ because outside 
the branching locus $Y$ may be covered by analytically open subsets 
biholomorphic to open subsets on $X$.
\end{itemize}

\noindent Finally, subfoliations $\mathcal{G}$ of the numerically trivial 
foliation $\mathcal{F}$ of $\alpha$ are always numerically trivial since test 
forms for $\mathcal{G}$ whose support is not intersecting 
$\mathrm{Sing\ } \mathcal{F}$ are also test forms for $\mathcal{F}$.

\noindent We use these observations to replace $X$ by a desingularization 
$\pi : \widehat{X} \rightarrow X$ of the indeterminacy locus of the 
pseudo-effective fibration and $L$ by $\pi^\ast L$. Let $A$ be an ample 
divisor on $X$ and $k \gg 0$. Then the curves of type
\[ D_1 \cap \ldots \cap D_{n-2} \cap F,\ D_i \in |kA|, \]
where $F$ is a fiber of the pseudo-effective fibration and $n = \dim X$, form a
family $(C_t)$ which is generically connecting on every fiber $F$.

\noindent There exists a composition $\pi : \widehat{X} \rightarrow X$ of 
modifications and a finite covering such that the strict transforms 
$(\overline{C_s})_{s \in S}$ of a subfamily of $(C_t)$ are the fibers of an 
everywhere defined holomorphic map $f: \widehat{X} \rightarrow S$ whose induced
foliation is contained in the numerically trivial foliation $\widehat{F}$ of
$\pi^\ast L$.

\noindent The general fiber curve $\overline{C_s}$ of $f$ is smooth and does 
not intersect $\mathrm{Sing\ } \mathcal{F}$ because the codimension of 
$\mathrm{Sing\ }(\widehat{\mathcal{F}})$ is $\geq 2$. Hence there is an 
analytically open set $V \subset S$ such that all fiber curves over $V$ are 
smooth and do not intersect $\mathrm{Sing\ }(\widehat{\mathcal{F}})$. 

\noindent Now we choose an $(n-1,n-1)$ test form $u$ for 
$\widehat{\mathcal{F}}$ on $U := f^{-1}(V) \subset \widehat{X}$ on which we can
apply Fubini's theorem:
\begin{eqnarray*}
0 & = & \lim_{\epsilon \downarrow 0} 
        \sup_{T \in \pi^\ast\alpha[-\epsilon\widehat{\omega}]} 
        \int_U (T_{ac} + \epsilon\widehat{\omega}) \wedge u = 
        \lim_{\epsilon \downarrow 0} \sup_T 
        \int_V \left( \int_{\overline{C_s}} T_{ac} + \epsilon\widehat{\omega}
               \right) \omega_V^{n-1} \\
  & = & \lim_{\epsilon \downarrow 0} \sup_T \int_V \left(
        \int_{\overline{C_s}} T - \sum_i \nu(T,E_i)[E_i] + 
                              \epsilon\widehat{\omega}
        \right) \omega_V^{n-1} \\
  & = & \deg \pi \cdot \lim_{\epsilon \downarrow 0} \sup_T \int_V \left(
        L.C_t - \sum_i \nu(T,E_i)E_i.C_t
        \right) \omega_V^{n-1} \\
  & \stackrel{(\ast)}{=} & \mathrm{Vol}(V) \cdot (L.C_t) -
    \lim_{\epsilon \downarrow 0} \sum_i \inf_T \nu(T,E_i)(E_i.C_t) =
    (L - N(L)).C_t.
\end{eqnarray*}

\noindent Here $\widehat{\omega}$ is a K\"ahler form on $\widehat{X}$ and 
$\omega_V$ a K\"ahler form on $V$. The equality $(\ast)$ is true because there
exists a sequence of currents $T_k \in \alpha[-\epsilon_k\omega]$, 
$\epsilon_k \stackrel{k \rightarrow \infty}{\longrightarrow} 0$, such that
\[ \nu(T_k, E) \rightarrow \nu(\alpha,E) \]
for all divisors $E$.

\noindent For the last claim suppose that for general points $x \in X$ there 
exist curves $C_x$ with $C_x.L = 0$ and $C_x$ is not contained in any fiber of 
$f$. Using the Chow variety as in \cite[2.1.2]{BCEKPRSW00} we conclude that 
there is a covering family $(C_t)_{t \in T}$ such that
\[ C_t . L = 0\ \mathrm{and\ } C_t \not\in \{\mathrm{fiber\ of\ } f\}. \]
There exists a composition $\pi : \widehat{X} \rightarrow X$ of 
modifications and a finite covering such that the strict transforms 
$(\overline{C_s})_{s \in S}$ of a subfamily of $(C_t)$ are the fibers of an 
everywhere defined holomorphic map $g: \widehat{X} \rightarrow S$, and still
\[ \overline{C_s} . \pi^\ast L = 0. \]
By the next proposition $g$ is numerically trivial hence by using the 
Key Lemma~\ref{KeyLemma} we see that $\pi \circ f$ is not the numerically 
trivial foliation of $\pi^\ast L$ contradicting the observations $(a)$ and 
$(b)$ from above.
\end{proof}

\begin{prop}
Let $f: X \rightarrow S$ be a fibration of curves $C_s$, $s \in S$, and $L$ a
pseudo-effective line bundle on $X$. If $C_s.L = 0$ the foliation induced by 
$f$ is numerically trivial w.r.t. $L$.
\end{prop}
\begin{proof}
Remind how we proved the Transversality Lemma~\ref{transvers-lem} by using the 
Cauchy-Schwarz and H\"older inequalities. This technique shows that is 
enough to check the numerical triviality condition $(NT)_u$ on those constant 
test forms $u$ which allow us to apply Fubini's theorem as in the proof before.
The proposition follows from 
\[ L.C_t \geq (L - N(L)).C_t. \]
\end{proof}

\begin{rem}
The proof of Prop.~\ref{pseff-fib-prop} shows that every covering family of 
curves $(C_t)$ such that the $C_t$'s
lie inside the fibers of the pseudo-effective fibration satisfies 
\[ C_t . (L - N(L)) = 0. \]
\end{rem}

\noindent This is not true for the partial nef reduction defined in 
\cite[\S 8]{BDPP04}. Nevertheless this reduction map is closely related to the 
pseudo-effective fibration. There are two differences: First the authors use as
numerical condition
\[ L . C_t = 0 \]
for the covering family $C_t$ defining the reduction map. We changed this 
condition to
\[ (L - N(L)).C_t = 0 \]
because we were not able to construct a defining family $(C_t)$ whose general 
member does not intersect the exceptional divisors in $N(L)$. Morally there 
should be such a family since exceptional divisors tend to be contractible.

\noindent Second the authors used Campana's reduction map instead of the 
generic reduction map. To get the same results as \cite{BDPP04} about
the Kodaira dimension of $L$ we have to apply Campana's reduction map on the 
fibers
of the pseudo-effective fibration of $L$ and use the properties of Boucksom's
divisorial Zariski decomposition.

\noindent In the nef case all these differences vanish:
\begin{prop}
If $L$ is nef the pseudo-effective fibration is the nef fibration. In 
particular the generic quotient of curves $C$ with $C.L = 0$ equals the Campana
quotient.
\end{prop}
\begin{proof}
For every curve $C$ in a fiber of the pseudo-effective fibration there is a
covering family of curves such that $C$ is the component of one of these 
curves. Nefness implies $L.C = 0$. The equality of generic quotient and Campana
quotient follows from \cite[Thm. 2.4]{BCEKPRSW00} which is of Key Lemma type.
\end{proof}

\bibliographystyle{alpha}

\newcommand{\etalchar}[1]{$^{#1}$}

\end{document}